\documentclass[final,onefignum,onetabnum]{siamart171218}

\usepackage{braket,amsfonts}

\usepackage{array}

\usepackage[caption=false]{subfig}
\usepackage{pgfplots}

\newsiamthm{claim}{Claim}
\newsiamremark{remark}{Remark}
\newsiamremark{hypothesis}{Hypothesis}
\crefname{hypothesis}{Hypothesis}{Hypotheses}

\usepackage{graphicx,epstopdf}

\Crefname{ALC@unique}{Line}{Lines}

\usepackage{amsopn}

\usepackage{xspace}
\usepackage{bold-extra}
\usepackage[most]{tcolorbox}

\colorlet{texcscolor}{blue!50!black}
\colorlet{texemcolor}{red!70!black}
\colorlet{texpreamble}{red!70!black}
\colorlet{codebackground}{black!25!white!25}

\usepackage[noprefix]{nomencl}  % for nomenclature

\usepackage{xifthen}
\usepackage{graphicx}
\usepackage{amsmath}
\usepackage{amsfonts}
\usepackage{amssymb}
\usepackage{amstext}
\usepackage{amsbsy}
\usepackage{epsfig}
\usepackage{epsf}
\usepackage{multirow}
\usepackage{float}
\usepackage{stmaryrd} 

\usepackage{xspace}
\usepackage{url}
\usepackage{verbatim}
\usepackage{ifthen}
\usepackage{rotating}
\usepackage{mathpazo}
\usepackage{bbding}
\usepackage{booktabs}

\usepackage{bbm}
\usepackage[toc,page]{appendix}
\usepackage[bbgreekl]{mathbbol}
\DeclareSymbolFontAlphabet{\mathbbl}{bbold}

\usepackage{algorithm,algpseudocode}
\usepackage{listings}
\usepackage[noprefix]{nomencl}  % for nomenclature

\newcommand{\grad}{\vec{\nabla}}

\newcommand{\norm}[2][{}]{\lVert#2\rVert_{#1}}
\newcommand{\vect}[1]{\mathbf{#1}}

\newcommand{\op}[1]{\mathbbl{#1}}
\newcommand{\kernel}[1]{\mathbbm{#1}}

\newcommand{\dx}[1][x]{\,d#1}
\newcommand{\domg}{\dx[\Omega]}

\newcommand{\sphere}{\ensuremath{\mathcal{S}}\xspace}
\newcommand*{\domain}{\ensuremath{\mathcal{D}}\xspace}
\newcommand*{\boundary}{\ensuremath{{\partial\domain}}\xspace}
\newcommand{\normal}{\ensuremath{\vec{n}}\xspace}
\newcommand{\bnormal}{\ensuremath{\normal_\mathrm{b}}\xspace}
\newcommand{\omgdnb}{\direction\cdot\bnormal}

\newcommand{\sn}[1][N]{\ensuremath{S_#1}\xspace}
\newcommand{\direction}{\ensuremath{\vec{\Omega}}\xspace}
\newcommand{\position}{\ensuremath{\vec{x}}\xspace}
\newcommand{\xslabel}[2][]{\ifthenelse{\isempty{#1}}{\mathrm{#2}}{\mathrm{#2},#1}}
\newcommand{\Sigt}[1][]{\ensuremath{\Sigma_{\xslabel[#1]{t}}}\xspace}
\newcommand{\Sigs}[1][]{\ensuremath{\Sigma_{\xslabel[#1]{s}}}\xspace}
\newcommand{\Sigf}[1][]{\ensuremath{\Sigma_{\xslabel[#1]{f}}}\xspace}
\newcommand{\unit}[1]{\,\mathrm{#1}}
\newcommand{\m}{\unit{m}\xspace}
\newcommand{\sfluxunit}{\,\ensuremath{\unit{cm}^{-2}\unit{s}^{-1}}\xspace}
\newcommand{\afluxunit}{\,\ensuremath{\unit{cm}^{-2}\unit{s}^{-1}\unit{st}^{-1}}\xspace}
\newcommand{\n}{d\xspace}
\newcommand{\M}{{N_d}\xspace}

\newcommand{\ma}{\mathcal{A}}
\newcommand{\mb}{\mathcal{B}}
\newcommand{\mf}{\mathcal{F}}
\newcommand{\mj}{\mathcal{J}}
\newcommand{\mpp}{\mathcal{P}}

\newcommand{\mg}{\mathcal{G}}

\newcommand{\vI}{\bold{I}}
\newcommand{\voL}{{\op{\bold{L}}}}
\newcommand{\oL}{{\op{{L}}}}
\newcommand{\voS}{{\op{\bold{S}}}}
\newcommand{\oS}{{\op{{S}}}}
\newcommand{\voF}{{\op{\bold{F}}}}
\newcommand{\oF}{{\op{{F}}}}
\newcommand{\voB}{{\op{\bold{B}}}}
\newcommand{\oB}{{\op{{B}}}}
\newcommand{\voI}{{\op{\bold{I}}}}

\newcommand{\ed}{\tilde{e}}

\lstdefinestyle{siamlatex}{%
  style=tcblatex,
  texcsstyle=*\color{texcscolor},
  texcsstyle=[2]\color{texemcolor},
  keywordstyle=[2]\color{texemcolor},
  moretexcs={cref,Cref,maketitle,mathcal,text,headers,email,url},
}

\tcbset{%
  colframe=black!75!white!75,
  coltitle=white,
  colback=codebackground, % bottom/left side
  colbacklower=white, % top/right side
  fonttitle=\bfseries,
  arc=0pt,outer arc=0pt,
  top=1pt,bottom=1pt,left=1mm,right=1mm,middle=1mm,boxsep=1mm,
  leftrule=0.3mm,rightrule=0.3mm,toprule=0.3mm,bottomrule=0.3mm,
  listing options={style=siamlatex}
}

\newtcblisting[use counter=example]{example}[2][]{%
  title={Example~\thetcbcounter: #2},#1}

\newtcbinputlisting[use counter=example]{\examplefile}[3][]{%
  title={Example~\thetcbcounter: #2},listing file={#3},#1}

\DeclareTotalTCBox{\code}{ v O{} }
{ %fontupper=\ttfamily\color{texemcolor},
  fontupper=\ttfamily\color{black},
  nobeforeafter,
  tcbox raise base,
  colback=codebackground,colframe=white,
  top=0pt,bottom=0pt,left=0mm,right=0mm,
  leftrule=0pt,rightrule=0pt,toprule=0mm,bottomrule=0mm,
  boxsep=0.5mm,
  #2}{#1}

\patchcmd\newpage{\vfil}{}{}{}
\begin{tcbverbatimwrite}{tmp_\jobname_header.tex}
\title{A highly parallel multilevel Newton-Krylov-Schwarz method with subspace-based coarsening and partition-based  balancing for the multigroup neutron transport equations on 3D unstructured meshes}

\author{Fande Kong\footnotemark[1]~\footnotemark[2] 
\and Yaqi Wang \footnotemark[3]
\and Derek R.~Gaston\footnotemark[2]
\and Cody J.~Permann\footnotemark[2]
\and Andrew E.~Slaughter\footnotemark[2]
\and Alexander D.~Lindsay\footnotemark[2]
\and Richard C.~Martineau\footnotemark[2]
}

\headers{A  parallel multilevel Newton-Krylov-Schwarz method}{Fande Kong}
\end{tcbverbatimwrite}
\input{tmp_\jobname_header.tex}

\begin{document}
\maketitle
\renewcommand{\thefootnote}{\fnsymbol{footnote}}
\footnotetext[1]{Corresponding Author: fande.kong@inl.gov} 
\footnotetext[2]{Computational Frameworks, Idaho National Laboratory, Idaho Falls, ID 83415} 
\footnotetext[3]{Nuclear Engineering Methods Development, Idaho National Laboratory, Idaho Falls, ID 83415} 
\renewcommand{\thefootnote}{\arabic{footnote}}

%% ------------------------------------------------------------------
%% ABSTRACT
%% ------------------------------------------------------------------
\begin{tcbverbatimwrite}{tmp_\jobname_abstract.tex}
\begin{abstract}
The multigroup neutron transport equations have been widely used to study the motion of  neutrons and their interactions with the background materials. Numerical simulation of the multigroup 
neutron transport equations is computationally challenging  because the equations is defined on a high dimensional phase space (1D in energy, 2D in angle, and 3D in spatial space), and furthermore, for realistic  applications, the computational spatial domain is complex  and the materials are heterogeneous.  The multilevel
domain decomposition methods is one of the most popular algorithms for solving the multigroup neutron transport equations, but the construction of coarse spaces is expensive and often not strongly scalable when the number of processor cores is large.   A  scalable algorithm
has to be designed in such a way that the compute time is almost halved without any comprise on the solution accuracy  when the number of processor cores is doubled.  In this paper, we study a highly parallel multilevel Newton-Krylov-Schwarz method equipped with  several  novel components, such as  subspace-based coarsening, partition-based balancing and hierarchical  mesh partitioning,  that  enable the overall simulation strongly scalable in terms of the compute time.   Compared with the traditional  coarsening method, the subspace-based coarsening algorithm  significantly reduces  the cost of the preconditioner setup that is often unscalable. In addition, the partition-based  balancing  strategy   enhances   the parallel efficiency of the overall solver by assigning a nearly-equal amount of work to each processor core. The hierarchical   mesh partitioning is able to generate a large number of subdomains and meanwhile minimizes  the off-node communication.  We numerically show that the proposed algorithm is scalable with more than 10,000 processor cores for a realistic  application  with a few billions unknowns on 3D unstructured meshes. 
\end{abstract}


\end{tcbverbatimwrite}
\input{tmp_\jobname_abstract.tex}
%% ------------------------------------------------------------------
%% END HEADER
%% ------------------------------------------------------------------

\section{Introduction}
The multigroup  neutron transport equations is employed to describe  the motion of neutrons and their interactions with the background materials \cite{lewis1984computational}.  The fundamental quantity  of interest  is the statistically averaged  neutron distribution, referred to  as ``flux",  in a high dimensional phase space (1D in energy, 2D in angle, 3D in space). We consider the time-independent version of the equations here so that the time dimension is not  taken into account.   The neutron flux is a scalar 
quantity   physically representing the total length traveled by all free neutrons per unit time and volume. For solving the neutron transport equations, some fundamental nuclear data (referred to as ``cross sections" ) describing the likelihood per unit path length of neutrons interacting with the background materials is  required. The cross sections depend on the energy and temperature of the background materials in a complicated manner \cite{lewis1984computational}.  The neutron transport equations can behavior  as  hyperbolic and elliptic forms under simple changes in cross sections (material properties) that may occur in realistic applications \cite{kaushik2009enabling}.  Because of the large dimensionality, the complicated solution behaviors, the complex computational domain and the heterogeneous  materials, the neutron simulations  are among the most memory and computation intensive in all of computational science.  Therefore, a scalable parallel solution approach  that takes advantages of modern supercomputers  plays  a
critical  role in the transport simulations.  In this paper, we propose a scalable parallel nonoverlapping  Newton-Krylov-Schwarz (NKS) method for the high-resolution simulation of the multigroup neutron transport equations.  The performance of NKS is almost  completely  determined  by the preconditioner. To achieve a high-performance neutron simulation, we develop a multilevel domain decomposition method with including  a novel   subspace-based  coarsening scheme, a  partition-based   balancing strategy and a hierarchical mesh partitioning approach. 

The development of efficient algorithms for the neutron simulations has been  an active research topic for  a couple of decades, and many solvers 
were studied.  Multilevel domain decomposition  and algebraic multigrid (AMG) methods are ones of  the most   popular algorithms for the numerical solution of the neutron transport equations. We briefly review the multilevel and AMG methods  here, and for other popular  methods such as the transport sweeps, interested readers 
are referred to \cite{lewis1984computational,  zeyao2004parallel}.  With the development  of supercomputers, the domains decomposition  methods 
become attractive because they are  naturally suitable for parallel computations. In \cite{van2011parafish}, the second-order even-parity form of the time-independent  Boltzaman transport equations is solved with  fGMRES preconditoned by an one-level overlapping domain decomposition method, where ILU together with CG is chosen as a local subdomain solver. The algorithm is numerically demonstrated  to scale up to a few hundreds processor cores, but the scalability drops significantly when using 1,000 processor cores.   In \cite{guerin2007domain}, a nonoverlapping domain decomposition with Robin interface conditions is studied for the simplified transport  approximation, and the parallel efficiency is reported using up to 25 processors cores.  In \cite{kaushik2009enabling}, the parallel computation is implemented using a space-angle-group decomposition  method, where  the ``within-group"  equation is solved  using a Richardson iteration.  A two-level overlapping Schwarz preconditioner is developed for the multigroup  neutron diffusions equations in \cite{kong2018fully}.

The algebraic multigrid (AMG) methods  can be implemented in space, angle or energy.   A spatial multigrid algorithm is presented for the isotropic neutron transport equations  with a simple 2D  geometry in \cite{chang2007spatial}, where the algorithm works well for homogeneous domains but the convergence need to be improved for heterogeneous  domains. 
In \cite{turcksin2012angular}, an angular multgrid method is used as a preconditioner for the GMRES method for the problems with highly forward-peaked scattering, and the method  is  numerically shown to be  more efficient than an analogous DSA-preconditioned (diffusion synthetic acceleration \cite{alcouffe1977diffusion}) Krylov subspace method.  A multigrid-in-energy preconditioner (MGE)  together with a Krylov subspace solver  is proposed in \cite{slaybaugh2013multigrid}. The MGE preconditioner reduces the number of Krylov iterations in both the fixed source and eigenvalue problems.

The multilevel and AMG methods are  used for a wide range of problems in the neutron transport, but the construction of coarse spaces is challenging and 
often unscalable when the number of processor cores becomes  large. To have scalable simulations, we take an attempt to address these issues using the NKS method equipped  with two important  ingredients,  subspace-based coarsening and partitioned-based  balancing. In this paper, the multigroup 
 transport equations is discretized   in space using the first order  continuous  finite element method and in angle using the discrete  ordinates  approach.  The resulting algebraic  system of equations   is  solved with Jacobian-free Newton-Krylov (PJFNK) \cite{knoll2004jacobian}, where the preconditioning matrix is formed with the streaming and collision operator. During each Newton iteration,  the Jacobian system is calculated by a Krylov subspace method such as GMRES preconditioned by a multilevel nonverlapping Schwarz. The coarse spaces can be constructed either  geometrically or algebraically.   In our previous works  \cite{kong2016highly, kong2017scalable, kong2018scalability}, some boundary preserving coarse spaces are constructed geometrically, and they are shown to work well for elasticity  problems and fluid-structure   interaction problems.  Unfortunately, the geometric coarsening   method is  unavailable  for the targeting application  since the computational domain, shown in Fig.~\ref{fig:mesh}, used in this work  includes many different regions that are meshed  using different element types. Instead,  an algebraic coarsening algorithm  is employed to   construct coarse spaces for the multilevel Schwarz method.  However,  if the traditional coarsening method is employed, the overall algorithm performance will be deteriorated and the strong scalability can not be  maintained. To overcome the difficulty, we introduce a novel subspace-based coarsening algorithm  that reduces the preconditionr setup time significantly  compared with the traditional coarsening method, which makes the overall algorithm scalable  with more than 10,000 processor cores.  In addition,  a partition-based balancing scheme  is included to enhance the parallel efficiency, and a hierarchical  mesh partitioning approach is studied to generate a large number of subdomains.

 The rest of this paper is organized as follows. In Section 2, the multigroup neutron transport equations and its spatial and angular discretizations  are described 
 in detail. And a highly parallel Newton-Krylov-Schwarz framework is presented in Section 3.  A novel subspace-based coarsening algorithm  is introduced, in Section 4,
 to construct coarse spaces for building an efficient Schwarz preconditioner. In Section 5, some numerical tests are carefully studied to demonstrate  the performance of the proposed algorithm.  A few remarks and conclusions are drawn in Section 6.

\section{Problem description}
In this section, we first describe the multigroup neutron transport equations in detail, and then present  the corresponding  spatial and angular discretizations. 
\subsection{Multigroup neutron transport equations}
\begin{figure}
\centering
    \includegraphics[width=0.38\linewidth]{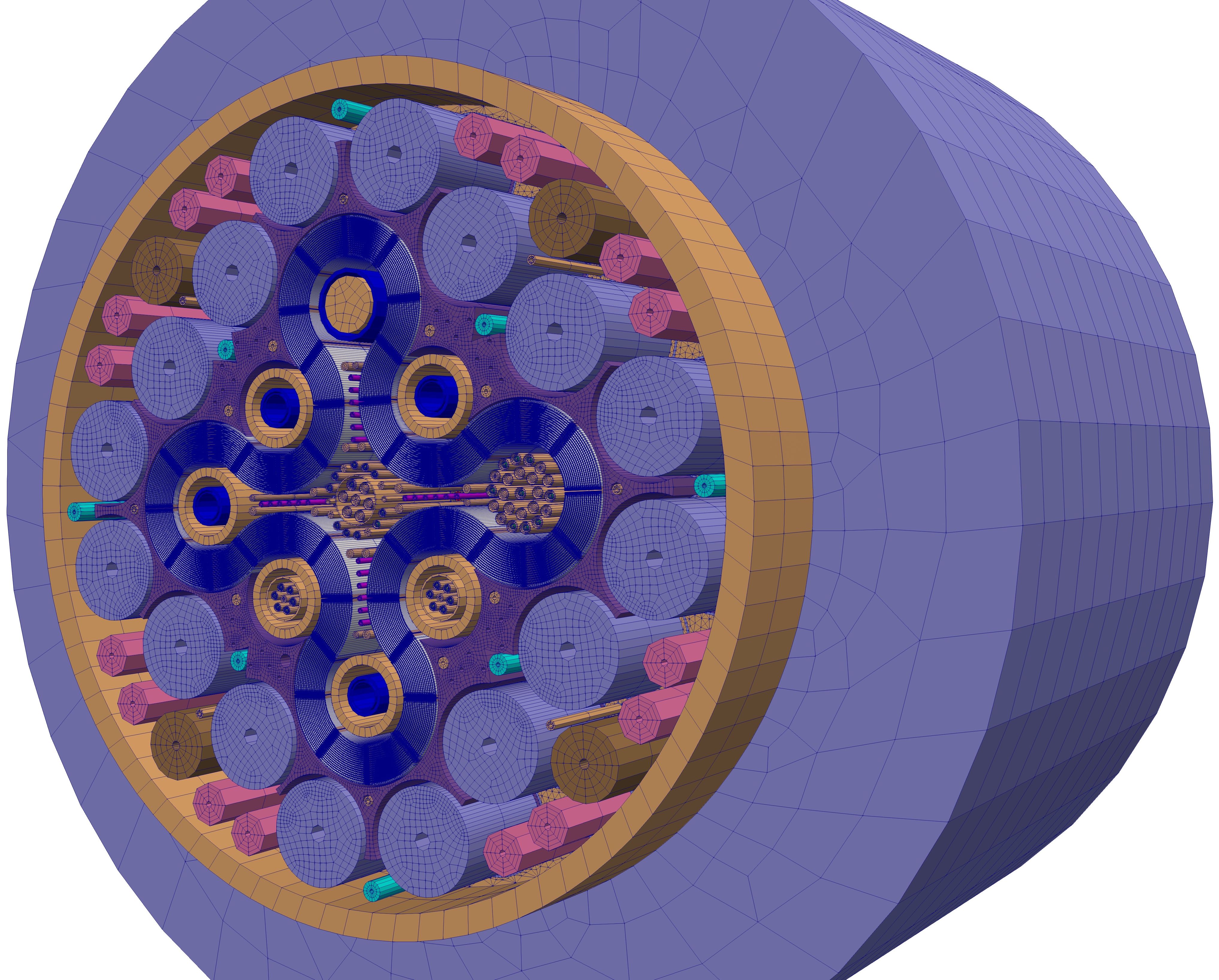} 
  \includegraphics[width=0.56\linewidth]{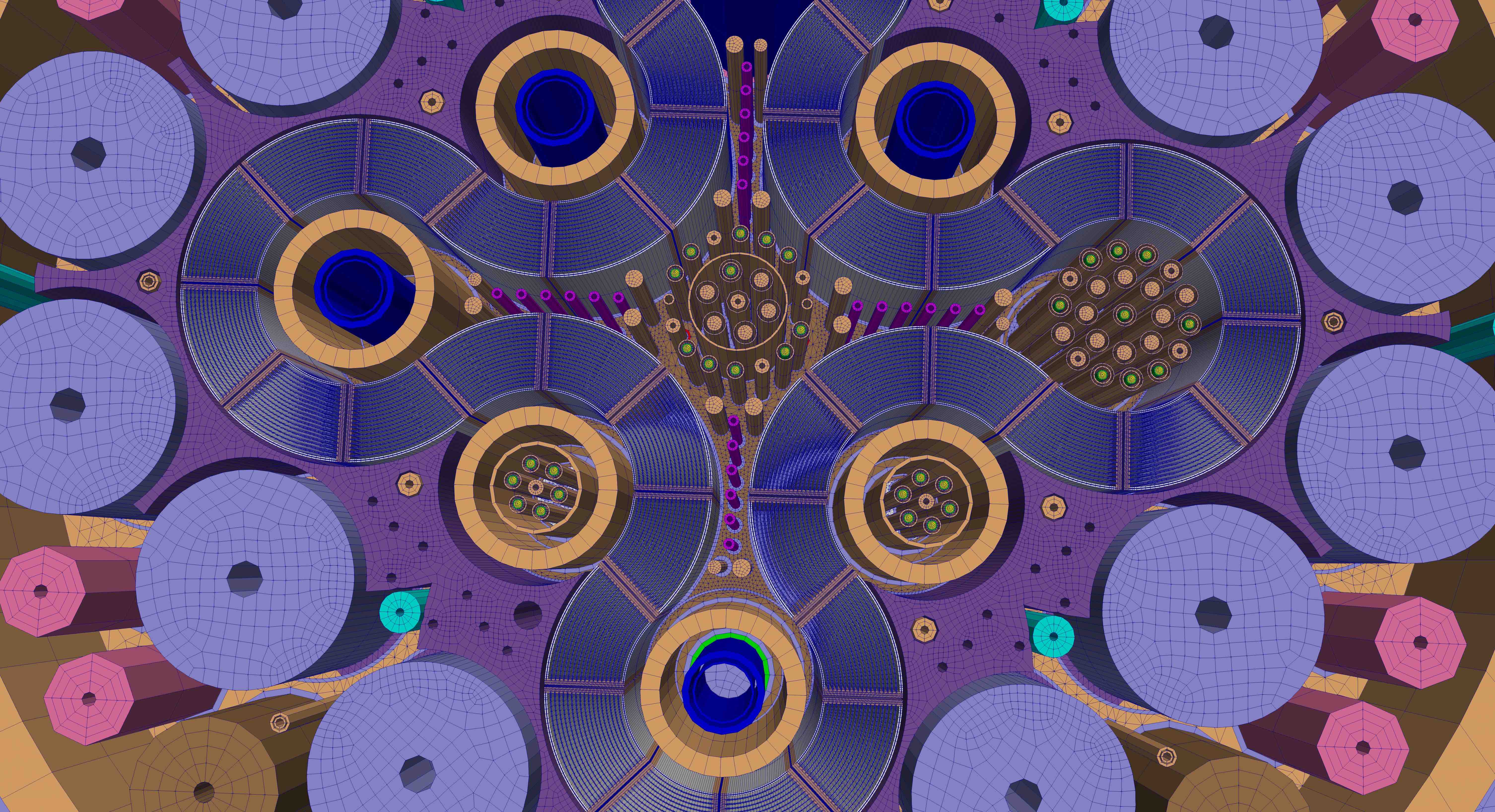} 
  \caption{3D unstructured  mesh. Right: zoom-in picture of the  top mesh. Different colors correspond to different materials.  \label{fig:mesh}}
\end{figure}

The fundamental  quantity  of interest, neutron angular   flux $\Psi_g$  $[\afluxunit]$,  is governed by the multigroup neutron transport  equations in $\domain \times \sphere$ as follows:
\begin{subequations}
\begin{equation}\label{eq:eigenvalue}
\begin{array}{llll}
\displaystyle
 \direction\cdot\grad\Psi_g +
 \Sigt[g]\Psi_g & = \displaystyle
 \sum_{g'=1}^G  \displaystyle \int_{\sphere} \Sigs[g'\rightarrow g]  f_{g'\rightarrow g}(\direction'\cdot\direction)\Psi_{g'}(\position, \direction')\domg'  \\
 & \displaystyle + \frac{1}{4\pi}\frac{\chi_g}{k}\sum_{g'=1}^G\nu\Sigf[g']\Phi_{g'}, \text{ in } \domain \times \sphere, \\
\end{array}
 \end{equation} 
 \begin{equation}\label{eq:eigenvalue_boundary}
 \Psi_g =
 \alpha_g^{\text{s}}\Psi_g(\direction_r) \displaystyle+
 \alpha_g^{\text{d}}\frac{ \displaystyle \int_{\direction'\cdot\bnormal > 0} \left| \direction'\cdot\bnormal \right|\Psi_g \domg'}
 { \displaystyle \int_{\direction'\cdot\bnormal > 0} \left| \direction'\cdot\bnormal \right|\domg'},
 \text{ on }\partial\domain:\omgdnb < 0,
 \end{equation} 
\end{subequations}
where  $g=1,\cdots,G$, and  $G$ is the number of energy groups. $\domain$ is a 3D spatial domain (e.g, shown in Fig.~\ref{fig:mesh})  and $\sphere$ is a 2D sphere.   $\position \in \domain$  is the independent spatial variable  $\left[ \unit{cm} \right]$, $\direction \in \sphere$  denotes the independent angular variable,  $\direction_r = \direction - 2(\omgdnb)\bnormal$,  $\partial \domain$ is the boundary of  $\domain$,  $\bnormal$ is the outward unit normal vector on the boundary,  $\Sigt[g]$ is the macroscopic total cross section $[ \unit{cm}^{-1} ]$,  $\Sigs[g'\rightarrow g]$ is the macroscopic scattering cross section from group $g'$ to group $g$  $[ \unit{cm}^{-1} ]$,  $\alpha^\text{s}_g$ is the specular reflectivity on $\partial \domain$,  $\alpha^\text{d}_g$ is the  diffusive reflectivity on $\partial \domain$, $k$ is  the eigenvalue (sometimes referred to as a multiplication factor),   $\chi_g$ is the prompt fission spectrum,  $\Sigf[g]$ is the macroscopic fission cross section $[\unit{cm}^{-1}]$,
and  $\nu$ is the averaged neutron emitted per fission.  $\Phi_{g}$ is the scalar flux $[\sfluxunit ]$ defined as $\Phi_{g}  \equiv \int_\sphere \Psi_g  \domg$, and $f_{g'\rightarrow g}$ is the scattering phase function.  In  \cref{eq:eigenvalue},  the first term  is the \emph{streaming term}, and the second  is  the \emph{collision} term.  The first term of~ \cref{eq:eigenvalue} on the right hand side is the \emph{scattering term}, which couples the angular fluxes of all directions and  energy  groups together. The second term of the right hand side of~\cref{eq:eigenvalue} is the fission term,  which also couples the angular fluxes of all directions and 
 energy  groups together.  For a more detailed description  on the neutron transport equations, please see~\cite{lewis1984computational, wang2018rattlesnake}.

For convenience, let us define some operators: 
%Streamingcollision
%\begin{equation}\label{eq:streamingcollision}
$$
\voL \vect{\Psi} \equiv
 \begin{bmatrix}
 \displaystyle
  \oL_1 \Psi_1, 
  \displaystyle
   \oL_2 \Psi_2,
   \displaystyle
   \cdots,
   \displaystyle
  \oL_G \Psi_G
\end{bmatrix}^T, ~ \oL_g \Psi_g \equiv   \grad \cdot \direction \Psi_g+ \Sigt[g] \Psi_g, 
$$
%\end{equation}
% Scattering
%\begin{equation}\label{eq:scattering}
$$
\voS \vect{\Psi} \equiv
 \begin{bmatrix}
\displaystyle
 \oS_1 {\Psi_1},
  \displaystyle
  \oS_2 {\Psi_2},
  \displaystyle
  \cdots,
  \displaystyle
\oS_G {\Psi_G}
 \end{bmatrix}^T,  ~ \oS_g {\Psi_g} \equiv  \sum_{g'=1}^G\int_\sphere \Sigs[g'\rightarrow g]f_{g'\rightarrow g}\Psi_{g'}\domg', 
 $$
%\end{equation}
% Fission
%\begin{equation}\label{eq:fission}
$$
\voF \vect{\Psi}  \equiv
 \begin{bmatrix}
 \displaystyle
\oF_1  \Psi_1 , 
  \displaystyle
\oF_2  \Psi_2  ,
  \cdots,
  \displaystyle
\oF_G  \Psi_G 
 \end{bmatrix}^T,  ~\oF_g  \Psi_g \equiv  \frac{1}{4\pi}\chi_{g}\sum_{g'=1}^G\nu\Sigf[g']\Phi_{g'}  .
 $$
%\end{equation}
Here $\voL$ is the streaming-collision operator,  $\voS$ is the scattering operator and $\voF$ is the fission  operator.  
% Boundary conditions
Similarly,  the operator for the boundary condition mapping from $\partial\domain\times\sphere^+_{\bnormal}$ to $\partial\domain\times\sphere^-_{\bnormal}$ is defined as 
%\hlabel{notation:boundary-reflecting}
%\begin{equation} \label{eq:boundary}
$$
 \voB \vect{\Psi} \equiv
  \begin{bmatrix}
   \displaystyle
 \oB_1  \Psi_1, 
     \displaystyle
 \oB_2  \Psi_2, 
     \displaystyle
    \cdots, 
     \displaystyle
     \oB_G  \Psi_G\\
  \end{bmatrix}^T, ~ \oB_g  \Psi_g \equiv   \alpha^\text{s}_g\Psi_g(\direction_r) +
    \alpha^\text{d}_g\frac{\int_{\direction'\cdot\bnormal > 0} \left| \direction'\cdot\bnormal \right|\Psi_g(\direction') \domg'}{\int_{\direction'\cdot\bnormal > 0} \left| \direction'\cdot\bnormal \right|\domg'}. 
$$
%\end{equation}
where $\sphere^\pm_{\bnormal} = \left\{\direction\in\sphere, \direction\cdot \bnormal  \gtrless 0\right\}$ is a half angular space defined with respect to the unit vector $\bnormal$.  Finally, 
 \cref{eq:eigenvalue} is rewritten as
 \begin{equation}\label{eq:operators}
\voL \vect{\Psi} = \voS \vect{\Psi} + \frac{1}{k}\voF \vect{\Psi},
\end{equation}
with the boundary condition $\vect{\Psi} =  \voB \vect{\Psi}$ corresponding to~\cref{eq:eigenvalue_boundary}.
\subsection{Spatial and angular discretizations} Before the weak  form of ~\cref{eq:operators} is presented, some notations are introduced.  An inner product is defined as
%\begin{equation}\label{eq:innerproduct}
$$
 \left( \vect{a}, \vect{b} \right)_{\domain\times\sphere} \equiv \sum_{g=1}^G\int_{\sphere}\domg\int_{\domain}\dx\, a_g(\position,\direction)b_g(\position, \direction),
$$
%\end{equation}
where $\vect{a}$ and $\vect{b}$ are generic multigroup functions defined in $\domain\times\sphere$.
We drop the subscript $\domain\times\sphere$ for notation simplicity.
We also have a similar  definition for  the  boundary integral  as :
%\begin{equation}\label{eq:surfaceinnerproduct}
$$
 \left\langle \vect{a}, \vect{b} \right\rangle \equiv  \left\langle \vect{a}, \vect{b} \right\rangle^{+} +  \left\langle \vect{a}, \vect{b} \right\rangle^{-},~
 \left\langle \vect{a}, \vect{b} \right\rangle^\pm \equiv
\sum_{g=1}^G \oint_{\boundary}\dx\int_{\sphere^\pm_{\bnormal}}\domg\, \left|\omgdnb\right| a_g(\position,\direction)b_g(\position, \direction). \\
$$
%\end{equation}
%\begin{equation}\label{eq:surfaceinnerproduct2}
% \left\langle \vect{a}, \vect{b} \right\rangle \equiv
%\sum_{g=1}^G \oint_{\boundary}\dx\int_{\sphere}\domg\, \left|\omgdnb\right| a_g(\position,\direction)b_g(\position, \direction),
%\end{equation}
Following  a standard finite element technique, we multiply a test function $\mathbf{\Psi}^\ast$ with~\cref{eq:operators}, and then  integrate over the phase space, $\domain\times\sphere$,
\begin{equation}
 \left(\vect{\Psi}^\ast, \voL \vect{\Psi}\right) = \left(\vect{\Psi}^\ast, \voS \vect{\Psi}\right)+ \frac{1}{k} \left (\vect{\Psi}^\ast , \voF \vect{\Psi} \right ).
\end{equation}
% Apply the property with  integration by parts,
%\begin{align}
%  \left(\vect{\Psi}^\ast, \op{L}\vect{\Psi}\right) + \left\langle \vect{\Psi}^\ast, \vect{\Psi}\right\rangle^- = \left(\op{L}^\ast\vect{\Psi}^\ast,\vect{\Psi}\right) + \left\langle 
%\vect{\Psi}^\ast, \vect{\Psi} \right\rangle^+, \label{eq:divergence-theorem}
%\end{align}
%on the streaming term, we have 
%\begin{align*}
%\left(\voL^\ast\vect{\Psi}^\ast, \vect{\Psi}\right)
% +\left\langle \vect{\Psi}^\ast, \vect{\Psi} \right\rangle^+
% -\boxed{\left\langle \vect{\Psi}^\ast, \vect{\Psi} \right\rangle^-} =
%\left(\vect{\Psi}^\ast, \op{S}\vect{\Psi}\right)
% +  \left (\vect{\Psi}^\ast , \frac{1}{k}\op{F}\vect{\Psi} \right ). 
%\end{align*}
%And  we further substitute the boundary condition in the boxed term and obtain
After some manipulations, the weak form reads  as
\begin{equation}\label{eq:unstabilized-bilinear}
\left(\voL^\ast\vect{\Psi}^\ast, \vect{\Psi}\right)
 +\left\langle \vect{\Psi}^\ast, \vect{\Psi} \right\rangle^+
 -\left\langle \vect{\Psi}^\ast, \voB \vect{\Psi} \right\rangle^- =
 \left(\vect{\Psi}^\ast, \voS \vect{\Psi}\right)
  + \frac{1}{k}  \left (\vect{\Psi}^\ast , \voF \vect{\Psi} \right ), 
\end{equation}
where $\voL^\ast$ is the adjoint operator of $\voL$.
The  form~\cref{eq:unstabilized-bilinear} is usually unstable, and here  a stabilizing technique, SAAF (self-adjoint angular flux),  is included to remedy this issue.  In the SAAF method, the streaming-collision operator $\voL$ is split   into two parts (the streaming operator $\voL_1$ and the collision operator $\voL_2$), 
\begin{align}
 \voL \vect{\Psi} \equiv \voL_1\vect{\Psi} + \voL_2\vect{\Psi},
\end{align}
where
%\hlabel{notation:streaming-collision-split}
%\begin{equation*}
$$
\voL_1\vect{\Psi} \equiv
 \begin{bmatrix} 
 \oL_{1,1} \Psi_1 ,
  \oL_{1,2} \Psi_2, 
  \cdots,
  \oL_{1,G} \Psi_G
 \end{bmatrix}^T,  ~ \oL_{1,g} \Psi_g  \equiv  \direction\cdot\grad \Psi_g,
 $$
 $$
 \voL_2\vect{\Psi} \equiv
 \begin{bmatrix}  
 \oL_{2, 1}  \Psi_1,
 \oL_{2, 2}  \Psi_2,
\cdots,
 \oL_{2, G}  \Psi_G
 \end{bmatrix}^T, ~ \oL_{2,g}  \Psi_g  \equiv   {\Sigt[g]}{\Psi_g}.
 $$ 
%\end{equation*}
The ``inverse" of $\voL_2$ is further defined  as 
%\hlabel{notation:inverse-collision}
\begin{equation*}
\voL ^{-1}_2\vect{\Psi} \equiv
 \begin{bmatrix}  
\oL_{2,1}^{-1} \Psi_1,
\oL_{2,2}^{-1} \Psi_2,
\cdots,
\oL_{2,G}^{-1} \Psi_G
 \end{bmatrix}^T, ~ \oL_{2,g}^{-1} \Psi_g  =  {\Psi_g}/{\Sigt[g]}. 
\end{equation*}
It is easy to verify that  $\voL^{-1}_2\voL _2=\voI$, $\voL^\ast_2 = \voL_2$, $\voL^\ast_1 = -\voL_1$ and $\voL_1\voL_2 = \voL_2\voL_1$.  With rearranging  ~\cref{eq:operators}, we have  
\begin{equation}\label{eq:afe}
 \vect{\Psi} = \voL^{-1}_2\left(\frac{1}{k}\voF \vect{\Psi} + \voS \vect{\Psi}  - \voL_1\vect{\Psi}\right),
\end{equation}
which is called  the \emph{angular flux equation} (AFE).
\nomenclature[abb]{AFE}{angular flux equation}%
We substitute~\cref{eq:afe} into the streaming kernel of~\cref{eq:unstabilized-bilinear},
\begin{align*}
 \boxed{\left(\voL_1^\ast\vect{\Psi}^\ast, \vect{\Psi}\right)}
 +\left(\voL_2^\ast\vect{\Psi}^\ast, \vect{\Psi}\right)
 +\left\langle \vect{\Psi}^\ast, \vect{\Psi} \right\rangle^+
 -\left\langle \vect{\Psi}^\ast, \voB \vect{\Psi} \right\rangle^-
 =\left(\vect{\Psi}^\ast, \voS \vect{\Psi}\right)
 +  \frac{1}{k} \left (\vect{\Psi}^\ast ,\voF \vect{\Psi} \right ),
\end{align*}
and obtain the following form after a few manipulations
\begin{align*}
 \boxed{\left(\voL_1\vect{\Psi}^\ast, \voL^{-1}_2\voL_1\vect{\Psi} \right)
 +\left(\voL_2^\ast\vect{\Psi}^\ast, \vect{\Psi}\right)
 +\left\langle \vect{\Psi}^\ast, \vect{\Psi} \right\rangle^+} 
 -\left\langle \vect{\Psi}^\ast, \voB \vect{\Psi} \right\rangle^-\\
 = \left(\voL^{-1}_2\voL \vect{\Psi}^\ast, \voS \vect{\Psi}\right)
+  \frac{1}{k} \left(\voL^{-1}_2\voL \vect{\Psi}^\ast, \voF \vect{\Psi} \right).
\end{align*}
We noticed that the boxed kernels are symmetric positive definite (SPD), and the calculation of the SPD system is possible using the multilevel method  equipped with algebraic coarse spaces. 
\nomenclature[abb]{SPD}{symmetric positive definite}%
SAAF is equivalent to SUPG (Streamline upwind/Petrov-Galerkin)~\cite{brooks1982streamline} with the inverse of group-wise total cross sections as the stabilization parameter.  Finally, we denote the weak form obtained using the SAAF method as
%\nomenclature[abb]{SUPG}{Streamline upwind/Petrov-Galerkin}%
\begin{equation}\label{eq:workform}
 \kernel{a}\left(\vect{\Psi}^\ast, \vect{\Psi}\right) = \frac{1}{k}  \kernel{f}\left(\vect{\Psi}^\ast, \vect{\Psi}\right),
\end{equation}
with 
\begin{align*}
  \kernel{a}\left(\vect{\Psi}^\ast, \vect{\Psi}\right) \equiv&
  \left(\voL_1\vect{\Psi}^\ast, \voL^{-1}_2\voL_1\vect{\Psi} \right)
 +\left(\voL_2\vect{\Psi}^\ast, \vect{\Psi}\right)
 +\left\langle \vect{\Psi}^\ast, \vect{\Psi} \right\rangle^+ \nonumber\\
& -\left(\voL^{-1}_2\voL\vect{\Psi}^\ast, \voS \vect{\Psi}\right)
 -\left\langle \vect{\Psi}^\ast, \voB \vect{\Psi} \right\rangle^-,\\
 \kernel{f}\left(\vect{\Psi}^\ast, \vect{\Psi}\right) \equiv&
 \left(\voL ^{-1}_2\voL \vect{\Psi}^\ast,  \voF \vect{\Psi} \right).
\end{align*}

The $\sn$ (discrete ordinates) method that can be thought of as a collocation method is considered for  the angular discretization.   Given an angular quadrature set $\left\{\direction_\n, w_\n, \n=1,\cdots,\M\right\}$ consisting of $\M$ directions $\direction_\n$ and weights $w_\n$, the multigroup transport equations is solved along these directions and all angular integrations in the kernels are numerically evaluated with the angular quadrature.  With the $\sn$ method, an integral of general functions over $\sphere$ is represented as a weighted summation, that is, 
$$
\int_{\sphere} \Psi_{g}\domg'  = \sum_{d=1}^{N_d} w_\n \Psi_{g,d}. 
$$
It is straightforward to apply the technique to~\cref{eq:workform}. Take the collision term as an example, we have 
\begin{equation}\label{eq:collision_sn}
 \left( \voL_2\vect{\Psi}^\ast,\vect{\Psi}\right) =\displaystyle  \sum_{g=1}^G \int_\sphere
\domg \left(\Sigt[g] \Psi_g^\ast, \Psi_g\right)_\domain =  \displaystyle \sum_{g=1}^G\sum_{\n=1}^\M w_\n \left(\Sigt[g] \Psi_{g,\n}^\ast, \Psi_{g,\n}\right)_\domain,
\end{equation}
where $(\cdot , \cdot)_{\domain}$ denotes that the integral is taken over $\domain$ only. For the spatial discretization, the first-order Lagrange finite element is applied to $(\cdot , \cdot)_{\domain}$.  For more details on the  angular and spatial discretization  of the neutron transport equations used in this work, please see \cite{wang2018rattlesnake}.
After the angular and spatial discretization, a large eigenvalue system with the dense  coupling block matrices in the  energy  and angle     is produced.  The potential
 dense matrix in energy  is generated because a high energy neutron can be scattered down to a low energy group  (down-scattering) and a low energy neutron can be also scattered up 
 to a high energy group (up-scattering). The equation is fully coupled in the angle.  We will introduce a scalable eigenvalue solver in next Section to 
 handle the large system of eigenvalue equations.

\section{Scalable parallel algorithm framework} In this Section, we describe the parallel  algorithm framework consisting of the Newton method for calculating  the nonlinear system of equations, the Krylov subspace method for solving the Jacobian  system and the Schwarz preconditioner for accelerating the linear solver.

 The corresponding agebraic system of equations for~\cref{eq:workform}  reads as 
\begin{equation}\label{eq:algebraic}
\ma \vect{\Psi} = \frac{1}{k}\mb \vect{\Psi}, 
\end{equation}
where  $\vect{\Psi}$ is also used to represent the solution vector that corresponds to the nodal values of the neutron flux at the mesh vertices, $\ma$ is the corresponding  matrix of  
$\kernel{a}$, and $\mb$ is the corresponding matrix of  $\kernel{f}$. Note that the matrices $\ma$ and $\mb$ are not necessary  to be formed 
explicitly, and we will have a detailed discussion on this shortly.  The simplest algorithm for the eigenvalue calculation of~\cref{eq:algebraic} is the inverse power iteration,  shown 
in  Alg.~\ref{alg:inverse_power_simplified}, that  works well only when the ratio of the minimum  eigenvalue to the second smallest eigenvalue is sufficient small, but it converges slow or even fails to converge  when the ratio is close to ``1".
% Inverse power iteration
\begin{algorithm} 
  \caption{Inverse power iteration.  ``$\leftarrow$" represents that the corresponding vector  is scaled in place. $\max_{e}$ is the maximum number of inverse power iterations. $\text{tol}_{\vect{\Psi}}$ and $\text{tol}_{k}$ are relative tolerances for the eigenvalue and the eigenvector, respectively.   \label{alg:inverse_power_simplified}}
  \begin{algorithmic}[1]  
   \item Initialize  $\vect{\Psi}_0$
   \item Compute eigenvalue:  $k_0 =  \norm{\mb \vect{\Psi}_{0} }$
   \item Scale $ \displaystyle \mb \vect{\Psi}_{0} \leftarrow \frac{1}{k_0} \mb \vect{\Psi}_{0} $
   \For {$n=0,..., \max_{e}$}
    \State $ \ma \vect{\Psi}_{n+1} = \mb \vect{\Psi}_{n}$
    \State  $k_{n+1}  = \norm{\mb  \vect{\Psi}_{n+1}}$
    \State Scale  $\displaystyle \mb \vect{\Psi}_{n+1}  \leftarrow \frac{1} {k_{n+1}} \mb \vect{\Psi}_{n+1} $ 
    \If {$\frac{\norm{\vect{\Psi}_{n+1}-\vect{\Psi}_{n}}} {\norm{\vect{\Psi}_{n}}} < \text{tol}_{\vect{\Psi}}$  and $\frac{|k_{n+1}-k_{n}|}{|k_{n}|} < \text{tol}_{k}$ }
     \State Break 
    \EndIf
   \EndFor
   \item  Output $k_{n+1}$ and $\vect{\Psi}_{n+1}$
  \end{algorithmic}
\end{algorithm}
%%%%%%%
The difficulty is overcome by  a Newton method that  accelerates  the convergence.   To take the advantage of Newton,
lines 5 and 6 of Alg.~\ref{alg:inverse_power_simplified} are rewritten as follows:
% Nonlinear system of equations for eigenvalue problems
\begin{equation}\label{nonlinear_equation_newton}
 \mf(\vect{\Psi})= 
\begin{array}{llll}
\displaystyle  \ma \vect{\Psi} - \frac{1}{\norm{\mb \vect{\Psi}}} \mb \vect{\Psi}. \\
\end{array}
\end{equation}
And then an inexact Newton is  applied to~\cref{nonlinear_equation_newton}. More precisely, for a given $\vect{\Psi}_n$, the new solution is updated as follows:
\begin{equation}\label{alg:newton_update}
 \vect{\Psi}_{n+1} =  \vect{\Psi}_{n} +  \bar{\vect{\Psi}}_{n}. 
\end{equation}
Here  $\bar{\vect{\Psi}}_{n}$ is the Newton update direction obtained by solving the following Jacobian system of equations
\begin{equation}\label{alg:newton_jacobian}
\mj(\vect{\Psi}_{n}) \bar{\vect{\Psi}}_{n}  = -\mf(\vect{\Psi}_{n}),
\end{equation}
where $\mj(\vect{\Psi}_{n})$ is the Jacobian matrix at $\vect{\Psi}_{n}$, and $\mf(\vect{\Psi}_{n})$ is the nonlinear function residual evaluated at $\vect{\Psi}_{n}$.  To save the memory, $\mj(\vect{\Psi}_{n})$ 
is not excplitly formed, instead, it is carried out in a matrix-free manner. The corresponding Newton is referred  to as ``Jacobian-free Newton" method~\cite{knoll2004jacobian}. That is,  a matrix-vector product, $\mj(\vect{\Psi}_{n}) \bar{\vect{\Psi}}_{n}$,  is approximated by 
$$
\mj(\vect{\Psi}_{n}) \bar{\vect{\Psi}}_{n} = \frac{\mf(\vect{\Psi}_{n}+ \delta  \bar{\vect{\Psi}}_{n}  ) - \mf(\vect{\Psi}_{n})}{\delta},
$$
where $\delta$ is a small permutation that is a square root of the machine epsilon in this paper.  \cref{alg:newton_jacobian} 
is  solved using an iterative method such as GMRES \cite{saad1986gmres}, and  a preconditioner is required to construct a scalable and efficient parallel solver. Let us rewrite~\cref{alg:newton_jacobian} 
as a  preconditioned form 
\begin{equation}\label{alg:newton_preconditioner}
\mj \mpp^{-1} \mpp  \bar{\vect{\Psi}}  = -\mf,
\end{equation}
where $\mpp$ is the preconditioning matrix that is often an approximation to $\mj$, and $\mpp^{-1}$ is a preconditioning process. The Jacobian $\mj$ is carried out in a matrix-free manner  since it has  dense diagonal blocks  since  all groups and all directions are coupled through the fission term and the scattering term.  The preconditioning matrix  $\mpp$ is formed explicitly  by only  taking into consideration  the first three terms of~\cref{eq:workform} since they form a SPD matrix that can be calculated using the multilevel method with  algebraic coarse spaces.  In fact, the angular fluxes in the energy and angle  are independent in the first three  terms of~\cref{eq:workform}. If the variables were ordered group-by-group and direction-by-direction, $\mpp$ is written as 
\begin{equation}\label{alg:precondition_matrix}
 \mpp = 
 \left [
\begin{array}{llllll}
\mpp_{0,0} &&& \\
& \mpp_{1,1} && \\
& &\ddots & \\
& &&\mpp_{G,G}  \\
\end{array}
\right ],
\end{equation}
where $P_{g,g}$ is a block diagonal matrix for the $g$th energy  group expressed as 
\begin{equation}\label{alg:precondition_matrix_direction}
 \mpp_{g,g} = 
 \left [
\begin{array}{llllll}
\mpp^{(0,0)}_{g,g} &&& \\
& \mpp^{(1,1)}_{g,g} && \\
& &\ddots & \\
& &&\mpp^{(N_d,N_d)}_{g,g}  \\
\end{array}
\right ].
\end{equation}
Here $G$ is the number of energy  groups, and $N_d$ is the number of angular directions per  energy group. $\mpp_{g, g'}$ represents the coupling matrix between groups $g$ and $g'$.  If the scattering and the fission 
terms were taken into account, $\mpp$ would be a fully coupled matrix instead of the diagonal  matrix shown in~\cref{alg:precondition_matrix}, that is, $\mpp_{g, g'} \neq 0 \text{~when~} g \neq g'$.
$\mpp_{g, g}^{(d, d')}$ represents the coupling matrix between angular directions $d$ and $d'$ in the $g$th group.  Similarly, if the fission term and the scattering term were considered  in the preconditioning  matrix, 
$\mpp^{(d,d')}_{g,g}$ would be a fully coupled dense matrix, that is,  $\mpp^{(d,d')}_{g,g} \neq 0 \text{~when~} d \neq d'$.  For the given group $g$ and  direction $d$, $\mpp^{(d,d)}_{g,g}$ is a large  sparse matrix 
obtained from the spatial discretization. It is easy to note that the block structure in~\cref{alg:precondition_matrix} corresponds to the multigroup approximation, and that in~\cref{alg:precondition_matrix_direction} corresponds to  the angular discretization.

\begin{figure}
\centering
  \includegraphics[width=0.40\linewidth]{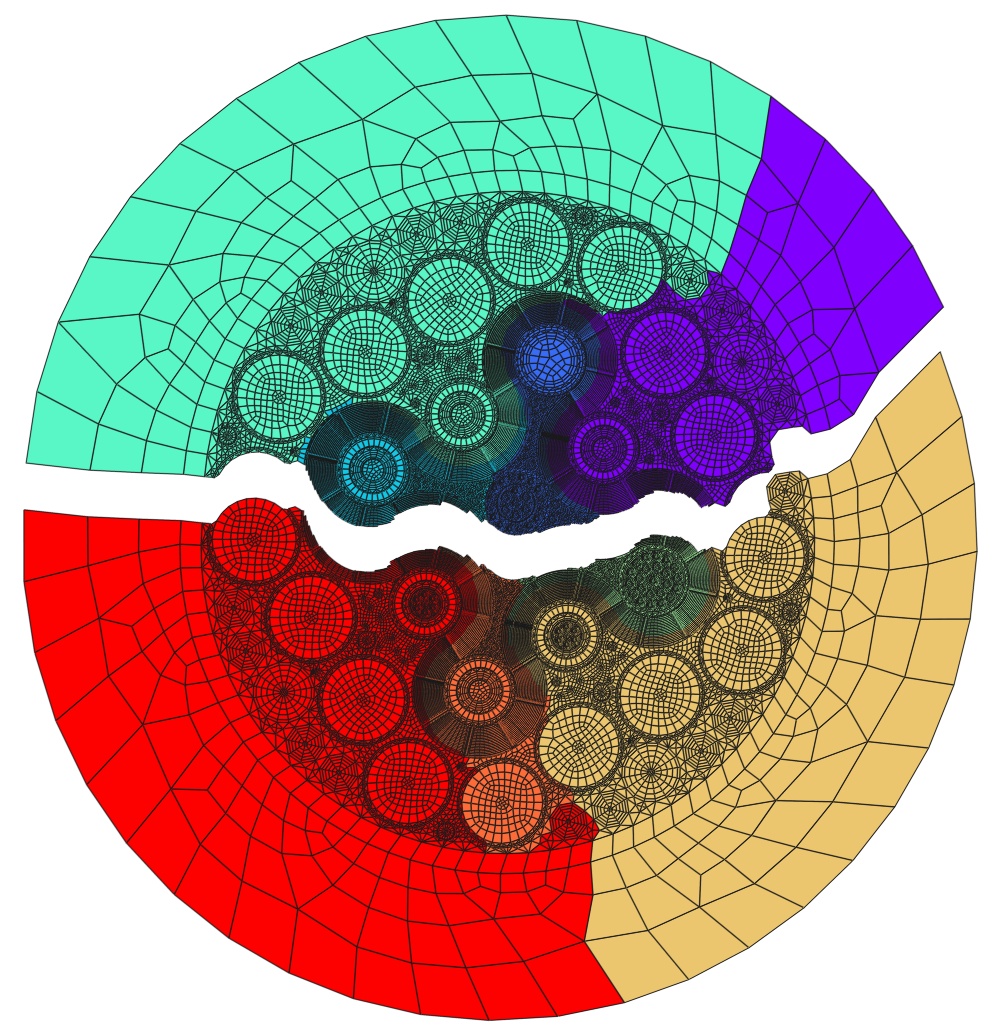} 
  \includegraphics[width=0.44\linewidth]{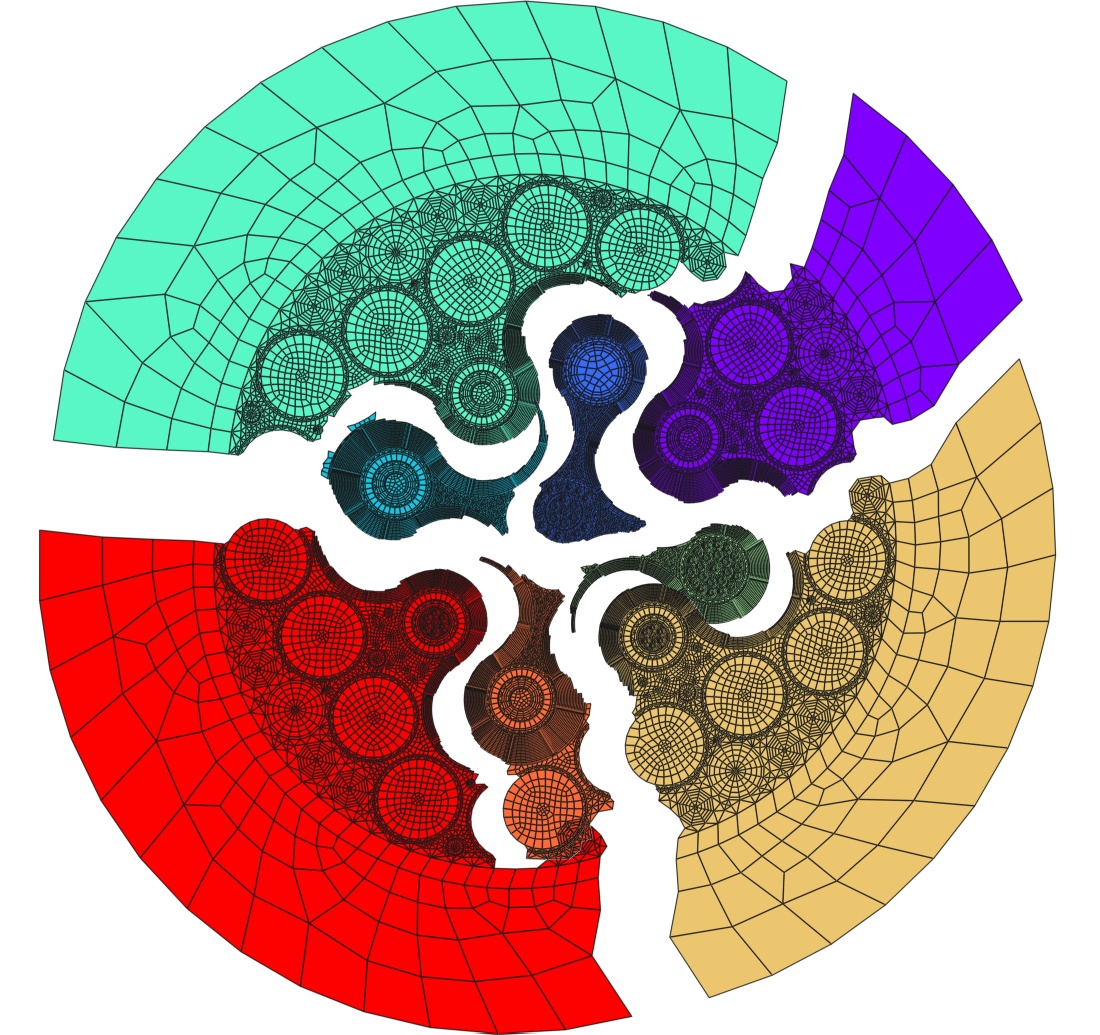} 
  \caption{Hierarchical partitioning. A mesh is paritioned into $2$ ``big" submeshes shown in the left, and then each ``big" submesh is further into 
  $4$ small submeshes. We finally have 8 submeshes. \label{fig:hierarchhical_partitioning_16}}
\end{figure}

Generally speaking, a preconditioning procedure is designed to find  the solution  of the following residual equations,
% Residual equation
\begin{equation}\label{alg:newton_preconditioner_process}
\mpp e = r,
\end{equation}
where $r$ is the residual vector from the outer solver (GMRES).  To carry out the simulation in parallel,  the mesh $\domain_h$,  corresponding to a triangulation  of $\domain$, is partitioned into $np$ ($np$ is the number of processor cores) submeshes $\domain_{h,i}, i=1, 2, .., np$.  This is accomplished by a hierarchical partitioning method since most 
existing partitioners such as ParMETIS \cite{karypis1997parmetis} do  not work well when the number of processor cores is close to or more than  10,000.  The basic idea of the hierarchical partitioning is to apply an existing partitioner such as ParMETIS or PT-Scotch \cite{chevalier2008pt} twice. The computational mesh $\domain_h$ is first partitioned 
into $np_1$ ``big" submeshes ($np_1$ often is the number of compute nodes), and each ``big" submesh is further divided into $np_2$ ($np_2$
is the number of processor cores per compute node) small submeshes. A 2D example with assuming that each compute node has $4$ processor cores is shown in Fig.~\ref{fig:hierarchhical_partitioning_16}, where the mesh is partitioned into $2$ ``big" submeshes,  and then 
each ``big" submesh is further divided into $4$ small submeshes, and finally we  have $8$ small submeshes in total.  Note that the hierarchical partitioning works for 3D meshes, and the 2D example is shown for the demonstration.
% partition figures
  Using the hierarchal partitioning method,  we are able to not only produce a large number of  submeshes, but also minimize the off-node communication since $np_2$ small submeshes on a compute node are physically connected and the communication between them is cheap.  Interested readers are referred to our previous works  \cite{kong2016highly,  kong2018general}  for more details of the hierarchical partitioning.  Let us denote the submatrix and the subvectors  associated with  a submesh $\domain_{h, i}$ as $\mpp_i$,  $e_i$ and $r_i$, respectively.  
We define a restriction operator, $R_i$, that restricts a global vector $r$ to  a nonoverlapping submesh, that is, $r_i$ = $R_i r$.   With those notations, the one-level nonoverlapping Schwarz  preconditioner  is expressed as 
\begin{equation}\label{alg:one_level}
\mpp^{-1}_{one} = \sum_{i=1}^{np} R_i^T \mpp_i^{-1} R_i, ~~ \mpp_i = R_i \mpp R_i^T
\end{equation}
where $\mpp_i^{-1}$ is a subdomain solver that is a successive over-relaxation (SOR) algorithm  in this paper.  The subdomain restriction $R_i$ does not extract any overlapping values. The overlapping version of~\cref{alg:one_level} has been successfully employed  in our previous works \cite{kong2016highly, kong2017scalable, kong2018efficient, kong2018simulation} for elasticity  equations, incompressible flows and fluid-structure interactions. Interested readers  are referred to \cite{smith2004domain, toselli2006domain} for more details  on the Schwraz  methods.
After many experiments, we find that the nonoverlapping Schwarz preconditioner is able to maintain a strong scalability for the targeting   applications when it is equipped with the subspace-based coarse spaces to be introduced shortly, and meanwhile the nonoverlapping Schwarz preconditioner uses less memory and  communication  compared with the overlapping version since no ghosting matrix entries  need to be stored and exchanged.   Coarse spaces need to be investigated  for  $\mpp^{-1}_{one} $ to form its multilevel version when the number of processor cores is large, the materials are heterogeneous  and the computational domain  is complex. Let us denote $L$  spaces  as $\domain_h^{(1)}, \domain_h^{(2)}, ..., \domain_h^{(L)}$, and the  associated operators as $\mpp^{(1)}, \mpp^{(2)}, ..., \mpp^{(L)}$. Here $\domain_h^{(1)} = \domain_h$ and $\mpp^{(1)} = \mpp$.   The interpolation operator from   $\domain_h^{(l+1)}$ to  $\domain_h^{(l)}$ is denoted as $I_{l+1}^l$ , and the corresponding restriction operator from $\domain^{(l)}_h$ to 
$\domain^{(l+1)}_h$ is $(I_{l+1}^{l})^T$.  A  multilevel additive Schwarz preconditioner (abbreviated as ``MASM")  is summarized  as  Alg.~\ref{MASM}.  
%% Multilevel method
 \begin{algorithm}
\caption{MASM($\mpp^{(l)}, e^{(l)}, r^{(l)}$)}\label{MASM}
\begin{algorithmic}[1]
\If {$l=L$} 
 \State Solve $\mpp^{(L)} e^{(L)} = r^{(L)}$ with a redundant direct solver  on each compute node
\Else 
 \State Pre-solve $\mpp^{(l)} e^{(l)} = r^{(l)}$ using an iterative solver preconditioned by  $(\mpp^{(l)}_{one})^{-1}$
 \State  Set $\bar{r}^{(l)} = r^{(l)} - \mpp^{(l)} e^{(l)}$
 \State  Apply the restriction: $ \bar{r}^{(l+1)} = (I_{l+1}^l)^T \bar{r}^{(l)}$
 \State  $z^{(l+1)}$ = MASM($\mpp^{(l+1)}, z^{(l+1)},  \bar{r}^{(l+1)}$)
 \State  Apply the interpolation: $z^{(l)} = I_{l+1}^l z^{(l+1)}$
 \State Correct the solution: $e^{(l)} = e^{(l)} + z^{(l)}$
 \State Post-solve $\mpp^{(l)} e^{(l)} = r^{(l)}$ using an iterative solver preconditioned by  $(\mpp^{(l)}_{one})^{-1}$
\EndIf
\item Return $e^{(l)}$ 
 \end{algorithmic}
\end{algorithm}
The fundamental  motivation behind  Alg.~\ref{MASM} is that the high frequency mode of the solution is efficiently  resolved using an iterative  method together with the preconditioner  $\mpp^{-1}_{one}$, and then the remaining low high frequencies  will be handled in the coarse levels.  The performance of Alg.~\ref{MASM} is largely affected by how to construct coarse spaces and their associated interpolations.  Generally speaking, there are 
two ways to construct a set of coarse spaces. The first approach  is to  geometrically coarsen the fine mesh $\domain_h$ to generate coarse meshes, which has been shown to be powerful 
in our previous works \cite{kong2016highly} for elasticity problems and  \cite{kong2016parallel, kong2017scalable, kong2018scalability} for fluid-structure interactions.  However, the geometry of  the targeting  application  is complex so that it is nontrivial to  setup a geometric mesh coarsening algorithm.  The second one is to construct coarse spaces without querying any mesh information, instead,  the coarse spaces and their interpolations are derived based on the matrix  information only. The second approach has been successfully applied for different applications   \cite{de2008distance, yang2002boomeramg, yang2006parallel}.    However,  it is  well known that the setup phase of the algebraic-version preconditioner  is not strongly scalable in terms of the compute time since the  matrix coarsening  and the interpolation construction are   expensive  \cite{yang2006parallel}. Fortunately, the overall algorithm can be still scalable if the preconditioner  setup phase accounts for a reasonably small portion of the total compute time.   We will introduce such a new subspace-based coarsening algorithm that 
the preconditioner setup time is significantly  reduced and the overall algorithm is able to  maintain a good scalability with more than 10,000 processor cores. We
will give a detailed description  of the proposed coarsening algorithm in next Section.

\section{Coarse spaces}  In this section, we discuss  a coarse space construction for Alg.~\ref{MASM}.  First, a matrix coarsening algorithm based on subspace is introduced, where
the ``grid" point selection  is accomplished using a submatrix instead of the global matrix.  A subinterpolation is constructed based on the splitting of the coarse points and the fine points, and the global interpolation is  built from the subinterpolation. 
\subsection{Matrix coarsening based on subspace}
According  to~\cref{alg:precondition_matrix} and~\cref{alg:precondition_matrix_direction}, it is easily found that $\mpp$ is a block diagonal  matrix and each block corresponds to the spatial discretization  of~\cref{eq:operators} for a given energy group and  angular direction.  Furthermore, there is no coupling between a block and the other blocks since we ignore the scattering and the fission terms in the preconditioning matrix.  The individual matrix blocks are similar to each other in the sense that they correspond  to the same continuous  operators and share the same computational mesh $\domain_h$. The differences between them come from different materials (i.e. cross sections) being used by different energy groups.  Our  motivation here is to 
coarsen a block of $\mpp$ instead of the entire matrix to generate  subinterpolations, and then the subinterpolations are expanded to covered the entire space by defining an expanding 
operator.  The benefit of this approach is  potentially save a lot of the setup time and also the memory usage since the coarsening phase operates on a much smaller data set.  Let us define a restriction $R_{s,i}$ that extracts the corresponding components from the entire vector $r$ for a given angular direction and  energy group  to form a subspace vector $r_{s, i}, i=1, 2, ..., G \times N_d$, that is,
$$
r_{s, i} = R_{s,i} r \equiv
\left ( 
\begin{array}{lll}
\vI & 0
\end{array}
\right ) 
\left (
\begin{array} {lll}
r_{s,i} \\
r/r_{s,i}
\end{array}
\right ),
$$
where ``$/$" denotes the components in $r$ but not in $r_{s,i}$.  The choice of energy  groups and angular directions is arbitrary in this paper, and we use the first energy group  and   angular direction, that is, $i=1$. Without any confusion,  we drop the second subscript of $r_{s, 1}$ and $R_{s,1}$,  and  denote them as $r_{s}$ and $R_{s}$, respectively, for the simplicity of notations. With these notations, a subspace preconditioning matrix (for the first energy group and  angular direction) is formed  as
%% sub matrix
\begin{equation}\label{eq:subspace_operator}
\mpp_s = R_s \mpp R_s^T. 
\end{equation}
%% sub matrix
Here $\mpp_s$ can be coarsened using one of the existing matrix coarsening algorithms.   We use a hybrid method of the Ruge-St$\ddot{u}$ben (RS) coarsening \cite{stuben2000algebraic}   and  the Cleary-Luby-Jones-Plassman (CLJP)
coarsening method \cite{yang2002boomeramg}.  For completeness, we briefly describe  these methods here, and interested readers are referred to \cite{luby1986simple, ruge1987algebraic, yang2006parallel} for more details. Before starting a coarsening process, a ``strength" matrix (graph), $\mg = (V, E)$,  need to be constructed from $\mpp_s$ since 
not all coefficients are equally important to determine the coarse spaces (grids) and we should consider the important coefficients only. Here $V$ is a set of all points  in $\mpp_s$, that is, $V=\{v_i\}$, and the size of $V$ is the number of rows of $\mpp_s$.  $E$ is a set of the corresponding edges, that is, $E=\{\ed_{ij}\}$.  An edge $\ed_{ij}$ is formed  when $v_i$ strongly depends on  $v_j$ or $v_j$ strongly influences   $v_i$ according to 
the following formula 
\begin{equation}\label{strong_dependence}
-p_{ij} \geq   \theta   \max_{k \neq i }(-p_{ik}),
\end{equation}
where $p_{ij}$ is an entity   of $\mpp_s$, and $\theta$ is the strength threshold that sometimes has an important  impact on the overall algorithm performance because it changes   the matrix complexities,  the stencil sizes,  and  the solver convergence rate. A coarsening algorithm tries  to split  $V$ into either  coarse points (C-point), denoted as $C$, which will be taken into the next level,  or fine points (F-points), denoted as $F$, which will be interpolated by C-points.   

The RS coarsening algorithm (also referred to as ``classical" coarsening in some literatures)  has two targets: 
\begin{enumerate}
 \item[A1]  For each point $v_j$ that strongly influences an F-point $v_i$,  $v_j$ is either a $C$-point  or it shares  a common C-point $v_k$ with $v_i$ 
 \item[A2]  $C$ should  be a maximal independent set
\end{enumerate}
``A1" is designed to insure the quality of interpolation, while ``A2" controls  the size of the coarse 
space  and  the complexity of the operator. In practice,  it is hardly to satisfy both  conditions at the same time. The RS coarsening  tries to meet A1 while uses A2 as a guideline and it  is carried out in two passes. In the first pass, each point $v_i$ is assigned by  a measure  $\m_i$ that equals  the number of the points strongly influenced by $v_i$, and  the point  with the maximum measure is  selected as 
C-point, $v_c$.  All the points  strongly influenced by $v_c$ are chosen as  new F-points, $\{v_f\}$. For each unmarked point that strongly influences any point in $\{v_f\}$, its measure is increased by the number of F-points  it influences. This procedure is repeated until all points are chosen as either C-points or F-points.  In the second pass, the algorithm checks every strong F-F connection  if two F-points  have  a common C-point.   If there is no a common C-point, and then one of the two F-points is chosen as a C-point.  The approach is summarized in Alg.~\ref{alg:RS_subspace}.  It is easily seen that the RS algorithm is inherently sequential. 
 \begin{algorithm}
\caption{Subspace based RS coarsening}\label{alg:RS_subspace}
\begin{algorithmic}[1]
\item Input: $\mpp$  \Comment{Submatrix extraction}
\item Extract a submatrix: $\mpp_s = R_s \mpp R_s^T$
\item Construct a strength matrix of $\mpp_s$, $\mg=(V, E)$, according to  ~\cref{strong_dependence}
\item Compute measures  $\{m_i\}$ for all points in $V=\{v_i\}$
\item Set $C=\emptyset$, $F=\emptyset$
\While {$V \neq \emptyset$}  \Comment{Pass 1}
\State Find a point $v \in V$ that has the maximum measure
\State $C=C+v$
\State Find the neighbors of $v$ (denoted as $V_n \subset V$) that strongly depend  on $v$
\State  $V=V-v$
\State $F=F+V_n$
\For {$v_n \in V_n$}
 \State Find the neighbors of $v_n$ (denoted as $V_{nn} \subset V $ ) that strongly influence $v_n$
  \For {$v_{nn} \in V_{nn}$}
  \State $m_{nn}$ = $m_{nn} + 1$
  \EndFor
\EndFor
 \State $V=V-V_n$ 
\EndWhile
\For {$v_i \in F$} \Comment{Pass 2}
\State Find  the neighbors of $v_i$ (denoted as $F_n \subset F$) that strongly influence $v_i$
\For {$v_n \in F_n$}
\If {$v_n$ and $v_i$ do not share a common C-point }
\State $F=F-v_n$
\State $C=C+v_n$
\EndIf
\EndFor
\EndFor
\item Output: $C$, $F$
 \end{algorithmic}
\end{algorithm}
A completely parallel coarsening approach is suggested in \cite{cleary1998coarse, yang2002boomeramg}. It is based on a  parallel maximal independent set  (MIS)
algorithm as described in \cite{luby1986simple}, and is often denoted as ``CLJP" (Cleary-Luby-Jones-Plassman)  in other literatures. The CLJP coarsening algorithm  starts with adding a measure $m_i$ for each 
point $v_i \in V$  just like  the RS coarsening algorithm.  Each $m_i$ is added by a small random value between $0$ and $1$ so that the points are distinctive  even if  the original 
measures are the same. It is now possible to find a local maximum of all the point measures  independently in parallel. A point $v_i$ with the local maximal measure is selected as a C-point, and the measures of  the neighboring  points strongly influenced by $v_i$ are decreased by 1. Furthermore, for all the points $\{v_j\}$ that strongly depend  on $v_i$, remove their  connections  to $v_i$.  Examine all the points $\{v_k\}$ that depend on $v_j \in \{v_j\}$ whether or not  they also depend on $v_i$. If $v_i$ is a common C-point of  $v_k$ and $v_j$,   remove the connection from $v_k$ to $v_j$ and decrease the measure of $v_j$ by 1. If the measure of $v_j \in \{v_j\}$ is smaller than 1, it is chosen as a F-point. This procedure is repeated until all points are selected as either C-points or F-points. The algorithm is summarized in Alg.~\ref{alg:CIJP_subspace}.
% CIJP coarsening
 \begin{algorithm}
\caption{Subspace based CIJP coarsening}\label{alg:CIJP_subspace}
\begin{algorithmic}[1]
\item Input: $\mpp$  \Comment{Submatrix extraction}
\item Extract a submatrix: $\mpp_s = R_s \mpp R_s^T$
\item Construct a strength matrix of $\mpp_s$, $\mg=(V, E)$, according to  ~\cref{strong_dependence}
\item Compute measure $\{m_i\}$ for all point in  $V=\{v_i\}$
\item Set $C=\emptyset$, and $F=\emptyset$
\item Add a random between 0 and 1 to each $m_i \in \{m_i\}$ \Comment{Pass 1}
\While {$V \neq \emptyset$}  
\State Find a point $v \in V$ that has the local maximum measure
\State $C=C+v$
\State Find the neighbors of $v$ (denoted as $V_n \subset V$) that strongly depend on $v$
\For {$v_n \in V_n$}
 \State $m_n = m_n -1$
 \State Find the neighbors of $v_n$ (denoted as $V_{nn} \subset V $ ) that strongly depends on $v_n$
  \For {$v_{nn} \in V_{nn}$}
  \If{$v_{nn}$ also depends on $v$}
  \State $m_{nn} = m_{nn}-1$
  \EndIf
  \EndFor
\EndFor
\For {$v_n \in V_n$}
  \If{$m_n < 1$}
  \State $F=F+v_n$
  \State $V=V-v_n$
  \EndIf
\EndFor
\State  $V=V-v$
\EndWhile
\item Output: $C$, $F$
 \end{algorithmic}
\end{algorithm}
% CIJP coarsening
While this approach works well for many applications, another option that has been shown to be even better is the combination of the RS coarsening and the CIJP coarsening \cite{yang2002boomeramg}.   This coarsening starts with Alg.~\ref{alg:RS_subspace} of which the first pass is applied to  the local graph  independently in parallel.   The   interior C-points and F-points   generated  in Alg.~\ref{alg:RS_subspace}   are used as an initial for  Alg.~\ref{alg:CIJP_subspace}. The resulting coarsening, which satisfies A1, fills the boundaries with further C-points and possibly adds a few in the interior of the subdomains. For convenience, the algorithm is denoted as ``HCIJP" (hybrid CIJP coarsening, and it is also referred to as ``Falgout" in \cite{yang2002boomeramg}) and shown in Alg.~\ref{alg:HCIJP_subspace}.
% CIJP coarsening
 \begin{algorithm}
\caption{Subspace based HCIJP coarsening}\label{alg:HCIJP_subspace}
\begin{algorithmic}[1]
\item Input: $\mpp$  \Comment{Submatrix extraction}
\item Extract a submatrix: $\mpp_s = R_s \mpp R_s^T$
\item Construct a strength matrix of $\mpp_s$, $\mg=(V, E)$, according to  ~\cref{strong_dependence}
\item Compute measure $\{m_i\}$ for all point in  $V=\{v_i\}$
\item Set $C=\emptyset$, and $F=\emptyset$
\item Apply the first pass of Alg.~\ref{alg:RS_subspace} 
\item Apply the first pass of  Alg.~\ref{alg:CIJP_subspace}
\item Output: $C$, $F$
 \end{algorithmic}
\end{algorithm}
% CIJP coarsening
 While these approach work well  for many applications, they sometimes  lead 
 to  high complexities.  There are some options that can be used to resolve these issues. The first one is to loose A1 as: A F-point should strongly depends on at least one C-point.  This approach often decreases the complexity, but the complexity  can  be still high and require more memory than desired.  This is further improved  by  an aggressive coarsening algorithm that is most efficiently implemented by applying the coarsening algorithms  twice, 
The resulting aggressive coarsening algorithm is briefly described in 
 Alg.~\ref{alg:AHCIJP_subspace}.
 %Aggressive hybrid  CIJP coarsening
 \begin{algorithm}
\caption{Subspace based aggressive HCIJP coarsening}\label{alg:AHCIJP_subspace}
\begin{algorithmic}[1]
\item Input: $\mpp$  \Comment{Submatrix extraction}
\item Extract a submatrix: $\mpp_s = R_s \mpp R_s^T$
\item Construct a strength matrix of $\mpp_s$, $\mg=(V, E)$, according to  ~\cref{strong_dependence}
\item Compute measure $\{m_i\}$ for all point in  $V=\{v_i\}$
\item Set $C=\emptyset$, and $F=\emptyset$
\item Apply the first pass of Alg.~\ref{alg:HCIJP_subspace}  to $\mg$
\item Apply the first pass of Alg.~\ref{alg:HCIJP_subspace}  to $C$
\item Output: $C$, $F$
 \end{algorithmic}
\end{algorithm}
% Aggressive hybrid  CIJP coarsening

\subsection{Interpolation construction based on subspace} 
With a splitting  $(C,F)$,  we consider the construction of  interpolation.   For a given F-point $v_i$, its interpolation takes the form  as follows:
$$
e_i = \sum_{j \in  C_i} w_{ij} e_j,
$$
where $C_i$ is the coarse interpolatory   set of $v_i$, and $w_{ij}$ is an interpolation weight determing the contribution of $e_j$ to $e_i$.  We 
assume  that an algebraically smooth error corresponds to a small residual, that is, $\mpp_s e \approx 0$ when $e$ is algebraically smooth. Let $V_{n,i}$ be the neighboring 
points of $v_i$, which strongly or weakly influence $v_i$, and then the $i$th equation of $\mpp_s e \approx 0$ reads as
\begin{equation}\label{interpolation_equation}
p_{ii} e_i + \sum_{j\in V_{n,i}}  p_{ij} e_j \approx 0.
\end{equation}
Here $V_{n,i}$ comprises three sets: $C_i$ (coarse neighbors),  $F_i^w$ (weakly influencing neighbors)  and $F_i^s$ (strongly influencing neighbors). A ``classical" 
interpolation, as described in \cite{ruge1987algebraic}, is constructed as
\begin{equation}\label{classical_interpolation}
w_{ij} = - \frac{1}{p_{ii} + \sum_{k\in F^w_i p_{ik}}} (p_{ij} + \sum_{k\in F^s_i} \frac{p_{ik}p_{kj}}{\sum_{m \in C_i p_{km}}}). 
\end{equation}
Eq.~\cref{classical_interpolation} is easy to implement in parallel since it only involves immediate neighbors and  only requires one layer of the ghosting  points. This method is invalid if A1 is not met.   Another interpolation scheme, often referred as ``direction interpolation", which only  needs  immediate neighbors and  can be used when A1 is violated,  is expressed as 
\begin{equation}\label{direct_interpolation}
w_{ij} = - \frac{p_{ij}}{p_{ii}} (\frac{\sum_{k \in V_{n,i}} p_{ik}}{ \sum_{l \in C_i}p_{il}}).
\end{equation}
If  an aggressive coarsening scheme such as Alg.~\ref{alg:AHCIJP_subspace} is adopted, it is necessary to use a long range interpolation, such as a  ``multipass interpolation''  as described in \cite{stuben2000algebraic}, in order to achieve a reasonable convergence. The ``multipass" interpolation scheme starts with computing interpolating weights using~\cref{direct_interpolation} for the F-points immediately influenced by at least one C-point. In the second pass, for each F-point $v_i$ that have not been interpolated yet, find its neighboring  interpolated F-points $v_j$ and then replace $e_j$ with $\sum_{k \in C_j} w_{jk} e_k$ in Eq.~\cref{interpolation_equation}. A direct interpolation is then applied to the modified equation. We would like to refer  interested readers to \cite{de2008distance}
for more details on different interpolation approaches. 
  Let us denote the subinterpolation constructed using the submatrix $\mpp_s$ as $I_{s,2}^1$ from the second  level to first  level. And the full interpolation $I_{2}^{1}$ is expanded using  the subinterpolation as follows
\begin{equation}\label{full_interpolation_equation}
I_{2}^{1} = \sum_{i=1}^{G \times N_d}(R_i^{(1)})^T I_{s,2}^1 R_i^{(2)} 
\end{equation}
where $R_i^{(l)}$ is the restriction  operator defined on the $l$th level for the $i$th variable.  The full coarse operator is  computed using a Galerkin method 
\begin{equation}\label{full_coarse_operators}
\mpp^{(2)} = (I_{2}^{1})^T \mpp^{(1)} I_{2}^{1}.
\end{equation}
With the full coarse operator in~\cref{full_coarse_operators} and the full  interpolation in~\cref{full_interpolation_equation}, the corresponding version of Alg.~\ref{MASM} is  denoted  
as ``MASM$_{\text{sub}}$", while that equipped with the  traditional coarse operators and interpolations   is simply written as ``MASM".  Note that  the description in this Section focuses  on generating one interpolation and one coarse operator, and  it is straightforward to apply the idea to generate a sequence of coarse spaces. 

\section{Numerical results}
\begin{figure}
\centering
  \includegraphics[width=0.3\linewidth]{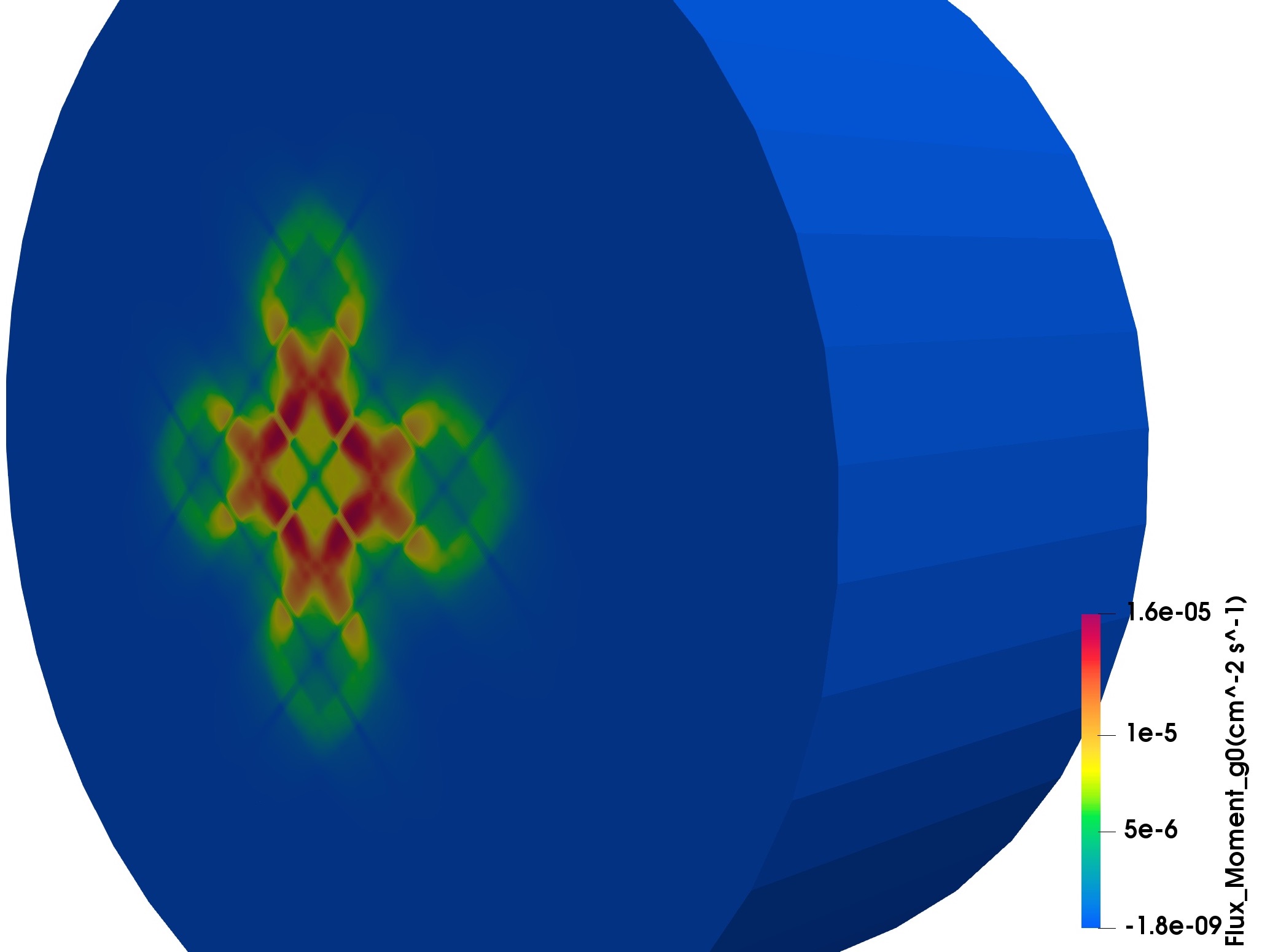} 
  \includegraphics[width=0.3\linewidth]{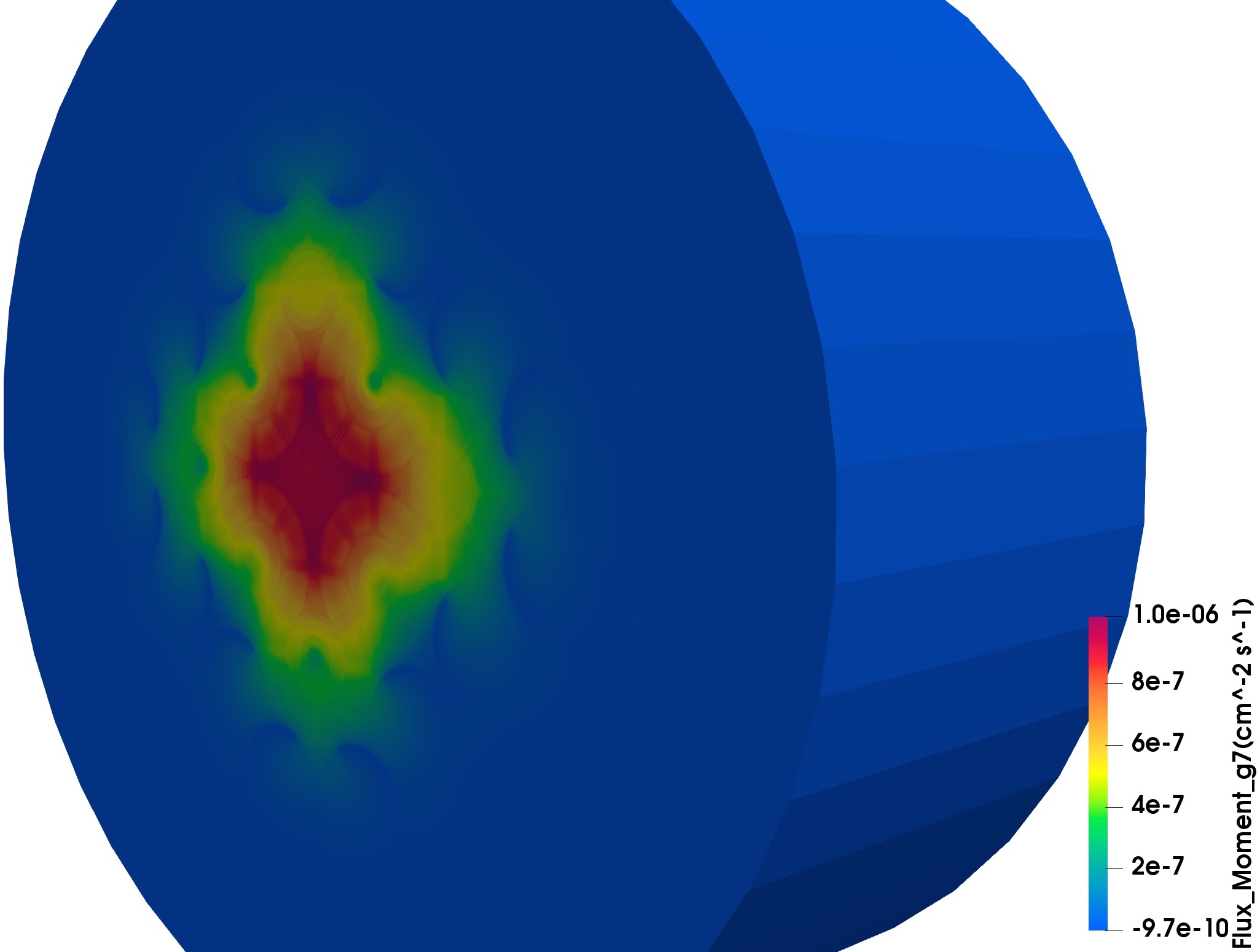} 
    \includegraphics[width=0.3\linewidth]{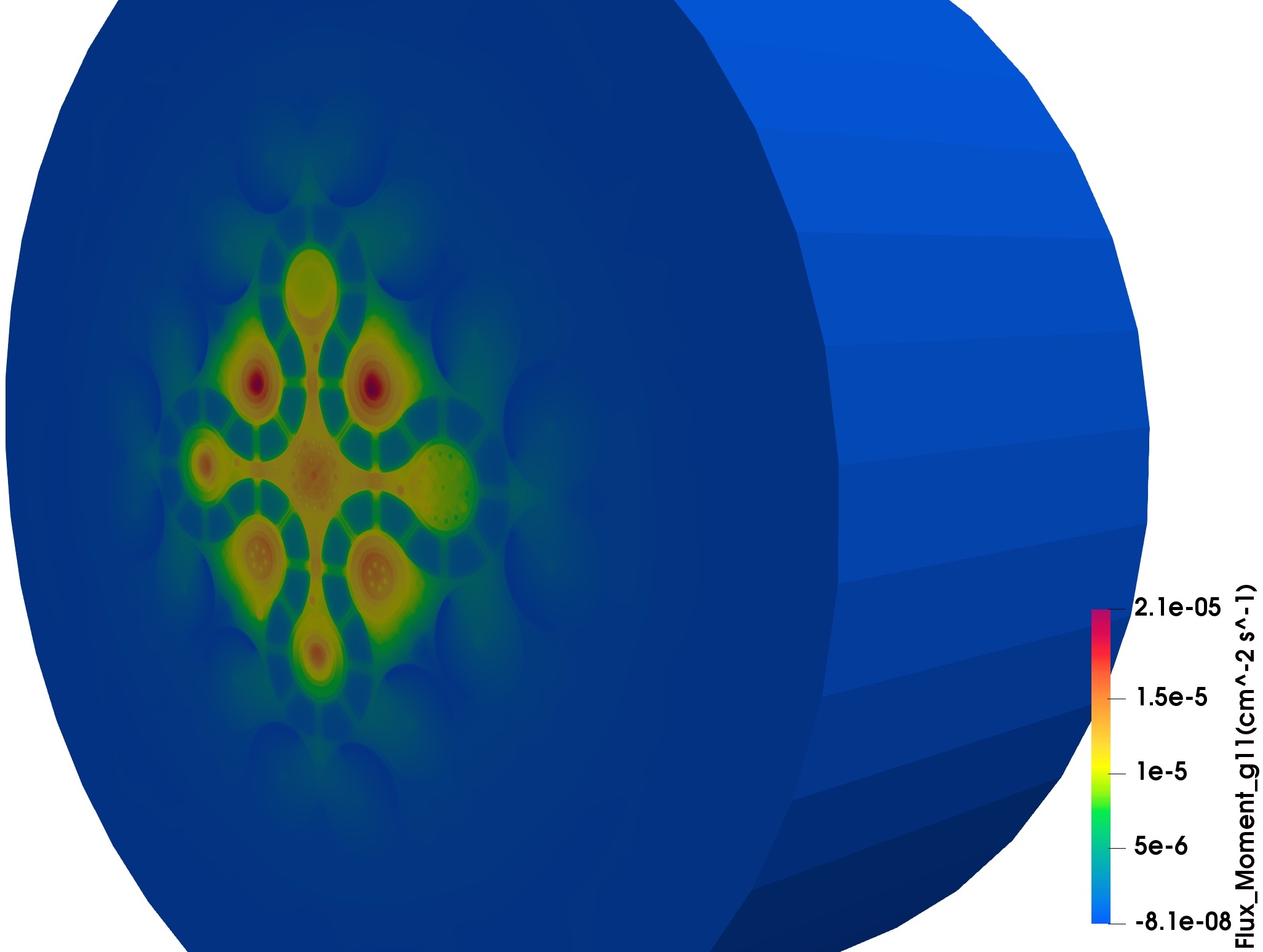} 
  \caption{Zero-order flux moments  for $1$st, $8$th and $12$th groups.  \label{fig:solution}}
\end{figure}

In this section, we report the algorithm's performance in terms of the compute time and the strong scalability  for the eigenvalue calculation of  the multigroup neutron transport equations for   a realistic application, namely  Advanced Test Reactor (ATR) that is located at the 
Reactor Technology Complex of the Idaho National Laboratory (INL) and is a 250-MW high flux test reactor. The ATR core, as shown in Fig~\ref{fig:mesh}, contains 40 fuel elements arranging  in a 
serpentine annulus between and around nine flux traps.  The  algorithms are implemented based on PETSc \cite{petsc-user-ref} and hypre \cite{hypre-web-page}. 
The numerical experiments are carried out on a supercomputer at INL, where each compute node has 36 processor cores (2.10 GHz per core) and the compute nodes  are connected by a FDR InfiniBand Network of 56 Gbit/s.  The problems are 
solved with an inexact Newton~\cref{alg:newton_update} together with GMRES preconditioned by Alg.~\ref{MASM}, where 4 iterations of the inverse power, as shown in Alg.~\ref{alg:inverse_power_simplified}, is used to generate an initial guess for Newton. In the Newton eigenvalue solver, a relative tolerance of $10^{-6}$
is enfored for the nonlinear solver, and an inexact linear solver with a relative tolerance of $0.5$ is adopted.  In the inverse power, one iteration of Newton together  with a linear solver with a relative tolerance of $10^{-2}$ is employed.  The eigenvalue functions for $1$st, $8$th and $12$th  are shown in Fig~\ref{fig:solution}. For convenience,  let us define some notations that will be used in the rest of discussions.  ``$np$" represents the number of processor cores,  ``NI" is the total number of Newton iterations, ``LI" denotes the total number of GMRES iterations, ``Newton" is the total compute time spent on the nonlinear solvers  and the inverse power iteration, ``LSolver" is the compute time on the linear solver, ``MF" is the compute time 
of  the matrix-free operations, ``PCSetup" is the compute time of  the preconditioner  setup, ``PCApply" is the compute time of the preconditioner apply,  ``EFF" is the parallel efficiency, and ``NR" is the ratio of the maximum number of  mesh nodes to the minimum  number of mesh nodes.  ``LSolver"  is part of ``Newton", and it  consists  of ``MF" and the preconditioner. The preconditioner time is split into ``PCSetup"  and ``PCApply".

\subsection{Comparison with traditional MASM}
\begin{figure}
\centering
  \includegraphics[width=0.8\linewidth]{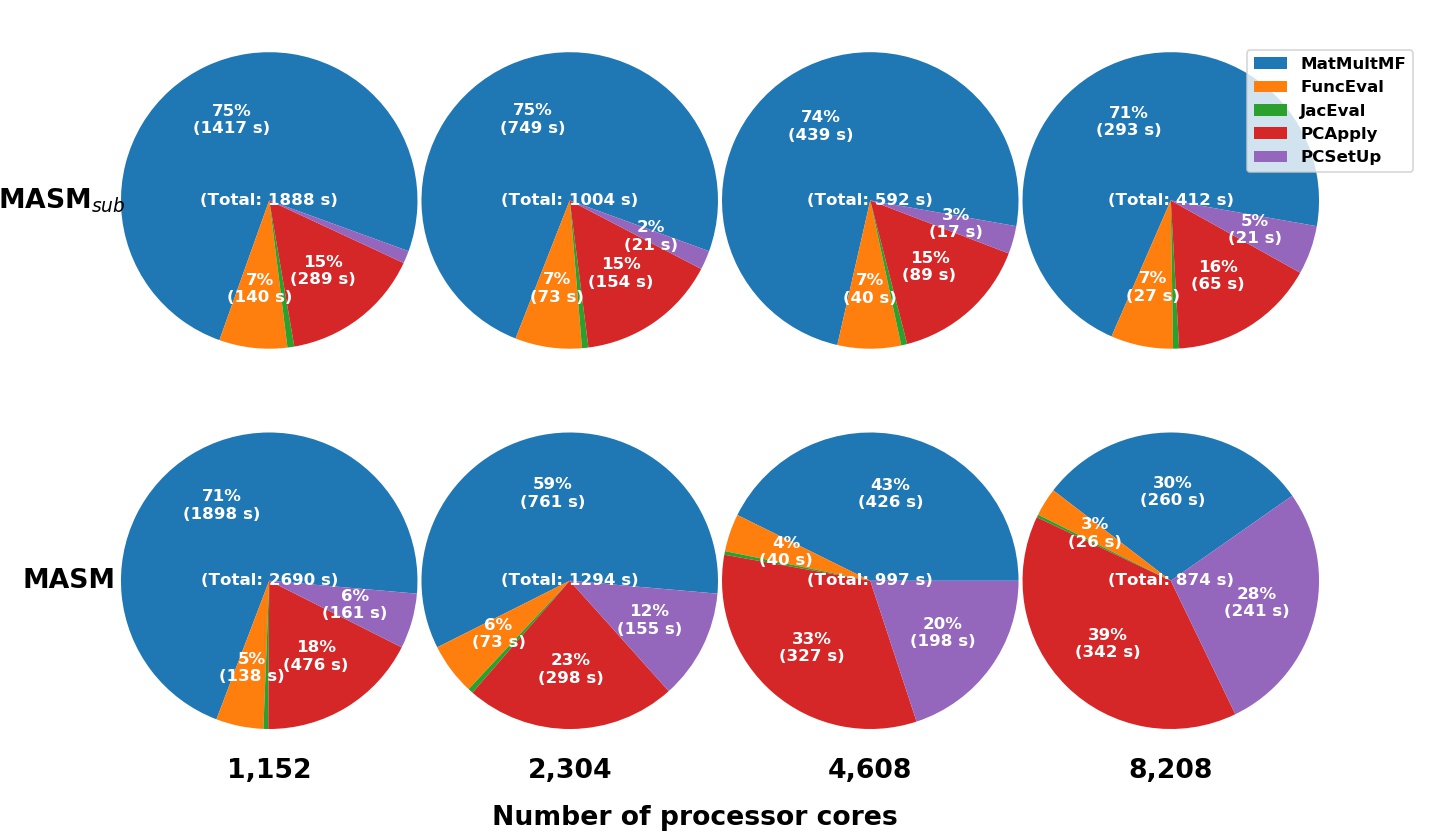} 
  \caption{Performance comparison with MASM. The figure shows the compute time spent on Jacobian evaluations (``JacEval"), function evaluations (``FuncEval"), matrix vector multiplications   via matrix-free (``MatMultMF"),  preconditioner setup (``PCSetUp") and preconditioner apply (``PCApply"), respectively.   \label{fig:pieplot_16_s2}}
\end{figure}

\begin{figure}
\centering
\includegraphics[width=0.32\linewidth]{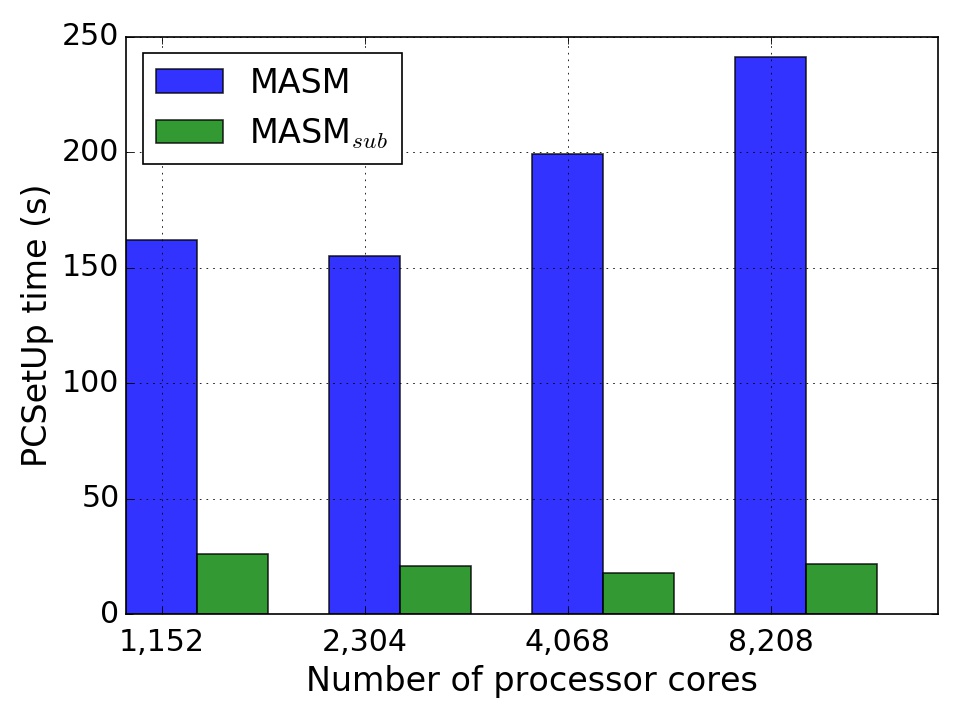} 
\includegraphics[width=0.32\linewidth]{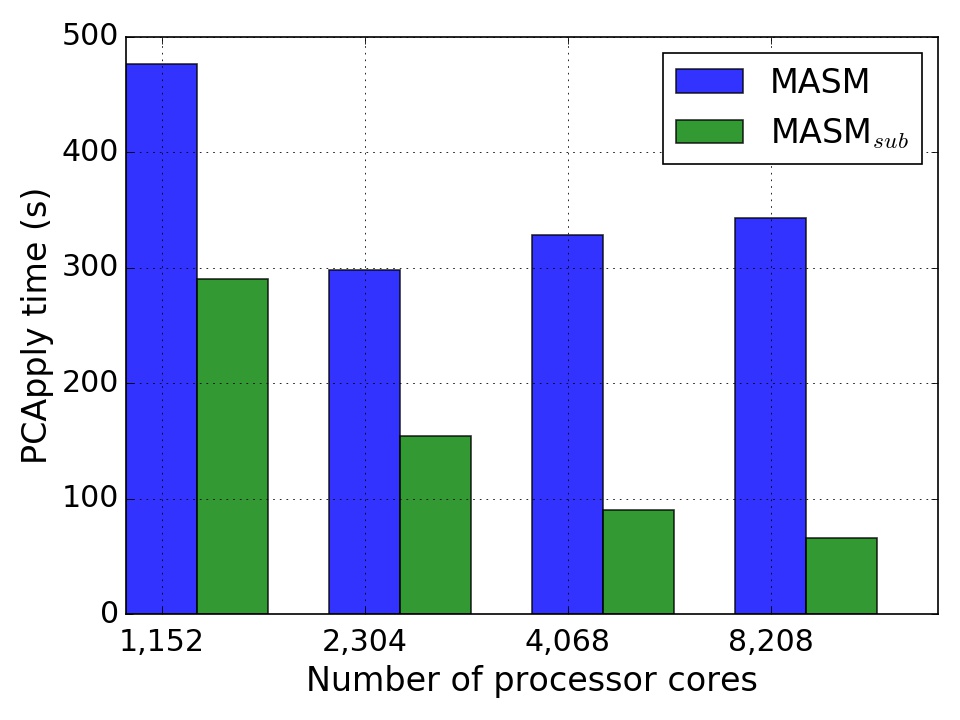} 
\includegraphics[width=0.32\linewidth]{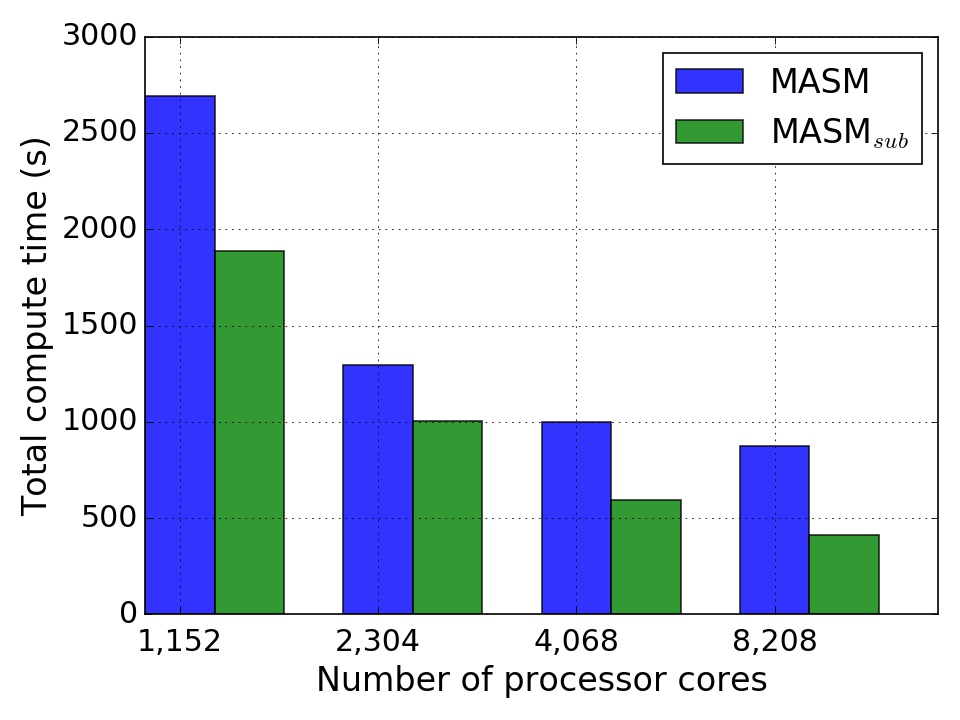}
\caption{The compute time comparison on different phases for the problem with 417,800,160 unknowns.  Left: preconditioner setup time; middle: preconditioner apply time; right: total compute time. \label{fig:barplot_16_s2}}
\end{figure}

We compare the proposed algorithm (denoted as ``MASM$_{\text{sub}}$") with the traditional MASM. We use a mesh with 4,207,728 elements and 4,352,085 nodes, where, at each node, there are $96$ unknowns consisting of $12$ energy groups and 8 angular directions. That is, the angle  is discretized by Level-Symmetric 2 with 8 angular directions.  The resulting system of nonlinear equations with 417,800,160 unknowns 
is solved using 1,152, 2,304, 4,608, and 8,208 processor cores, respectively. The performance comparison with the traditional MASM is summarized in Table~\ref{tab:amg_sub_coarse_s2} and Fig.~\ref{fig:pieplot_16_s2}.
% | grep "Linear |R| =" | wc -l 
% Performance comparison
\begin{table}
\scriptsize
\centering
\caption{Performance comparison with MASM for a problem with 417,800,160 unknowns. The resulting system of nonlinear equations with 417,800,160 unknowns 
is solved by an inexact Newton with MASM$_{\text{sub}}$ and  MASM  on 1,152, 2,304, 4,608, and 8,208 processor cores, respectively.\label{tab:amg_sub_coarse_s2}}
\begin{tabular}{c c c c c c  c c c c}
\toprule
$np$  &scheme& NI& LI& Newton& LSolver & MF   &  PCSetup & PCApply   & EFF \\
\midrule
1,152 & MASM$_{\text{sub}}$ &  13 & 191 & 1855 & 1701&1418& 26 &290 & 100\%\\
1,152 & MASM               &  13 & 251 & 2640 & 2486&1900& 162 &  476 & -- \\
\midrule
2,304 &  MASM$_{\text{sub}}$  & 13 & 193& 989& 908 & 749 &21 & 154 & 93\%\\
2,304 &  MASM  & 13 & 196& 1277& 1196 & 761 &155 & 298 & 73\%\\
\midrule
4,608 &  MASM$_{\text{sub}}$  & 13 & 202& 581& 535 & 440 &18 & 90 & 80\%\\
4,608 &  MASM  & 13 & 194& 985& 939 & 426 & 199 &  328 & 47\%\\
\midrule
8,208 &  MASM$_{\text{sub}}$  & 14 & 216& 404& 372 & 294 &22 & 66 & 64\%\\
8,208 &  MASM  & 13 & 192& 866& 835 & 261 &241 & 343 & 30\%\\
\bottomrule
\end{tabular}
\end{table}
% Performance comparison with MASM for S2
%%Discussion on the numerical  results
The nonlinear eigenvalue solver consists of Jacobian evaluation, function evaluation, matrix-vector multiplication, preconditioner setup and preconditioner apply,  and where most components except the  preconditioner setup  are mathematically scalable. As we discussed   earlier,    the  preconditioner setup including the matrix coarsening and the interpolation construction is challenging  to parallel, and its compute time sometimes  increases significantly when we increase the number of processor cores, which deteriorates  the overall algorithm.  From Fig.~\ref{fig:pieplot_16_s2},  we observed that the preconditioner setup  for the traditional MASM is not scalable, and the ratio of the  preconditioner setup time to the total compute time is increased significantly when we increase the number of processor cores. The ratio is only $6\%$ when the number of processor cores is $1,152$, but it jumps to $28\%$ when we use $8,208$ processor cores.  For the preconditioner setup time, it is $161~s$ at $1,152$ cores and increased to $261~s$ when  $8,208$ processor cores is used. In the traditional MASM, the precodnitioner setup not only is  unscalable, but also takes  a big chunk of the total compute time so that the overall algorithm performance is deteriorated and the parallel efficiency is reduced to  $30\%$ at $8,208$  processor cores. On the other hand,  the preconditioner  setup  of MASM$_{\text{sub}}$ performs better since it accounts for  only $3\%$ ($26~s$) of the total compute time at $1,152$ cores and it slightly increases  to $6\%$ ($22 ~s$) when we use $8,208$ processor cores. An interesting thing is that the preconditioner setup time of MASM$_{\text{sub}}$ does not increase much and stays  close to a constant. That makes the overall algorithm scale much better, and the parallel efficiency is about $64\%$ even when the number of processor cores is large, i.e., $8,208$.  The coarsening algorithm affects  not only the preconditioner setup time but also the preconditioner apply  time. In the traditional MASM, we observed that the preconditioner apply is not ideally scalable since while it accounts  for $16\%$ of the total compute time for $1,152$ processor cores, the ratio is increased to $39\%$ at $8,208$ cores.  The corresponding preconditioner apply  time is $476 ~s$ for $1,152$ processor cores, and it is decreased to $298 ~s$ by $37~s$ when we
double the number of processor cores.   Ideally, the preconditioner apply  time should be  reduced by  $50\%$  when the core count is doubled. The preconditioner apply  of MASM$_{\text{sub}}$ is scalable in the sense that the compute time is decreased from $289~s$ to $154~s$ by $47\%$ when we double the number of processor cores from 1,152 to 2,304, and it is further decreased to $89~s$ when we use 
4,608 processor cores. The traditional MASM does not preserve this property, and its preconditioner apply time is actually increased to $327~s$ from $298~s$ when the core count is doubled  from 2,304 to 4,608.  The coarsening algorithm based on subspace make MASM$_{\text{sub}}$ scalable for the ATR simulation while the traditional MASM does not perform well. At 8,208 core, MASM$_{\text{sub}}$ is twice faster than MASM.  These behaviors can be observed from Table 1 as well, where the number of Newton iterations is similar for both MASM and MASM$_{\text{sub}}$, and the GMRES iteration of MASM$_{\text{sub}}$ is slightly more than that of  MASM at 4,608 and 8,208 cores. The impact of the slight increase of GMRES iteration is negligible since  the preconditioner apply per iteration  is scalable for MASM$_{\text{sub}}$.   ``LSolver" accounts for the most of  the overall compute time, and the overall algorithm is scalable as  long as the linear solver performs well. The peformance of the linear solver is almost completely determined by the preconditioner since ``MF" is well-known to be scalable mathematically.  In summary, the eigenvalue solver together with MASM is not scalable, while the MASM$_{\text{sub}}$ equipped eigenvalue solver performs well since the setup phase of MASM$_{\text{sub}}$ is optimized and the apply phase of MASM$_{\text{sub}}$ scales well.  The same performance comparison is observed in Fig.~\ref{fig:barplot_16_s2} as well, where the preconditioner setup of MASM$_{\text{sub}}$ is almost $10$ times faster than MASM for all processor counts. The preconditioner apply for MASM$_{\text{sub}}$ is $2$ or $3$ times more efficient than MASM when the number of processor cores is small, and $5$ times faster at 8,208 cores. Due to these behaviors, the overall algorithm based on  MASM$_{\text{sub}}$ is much better than that with MASM. Note that the total compute time in Fig.~\ref{fig:pieplot_16_s2} is slightly  more than that in Table~\ref{tab:amg_sub_coarse_s2} since it is calculated  by summing  up all the individual  components  that have  some  overlap. 

\begin{figure}
\centering
  \includegraphics[width=0.8\linewidth]{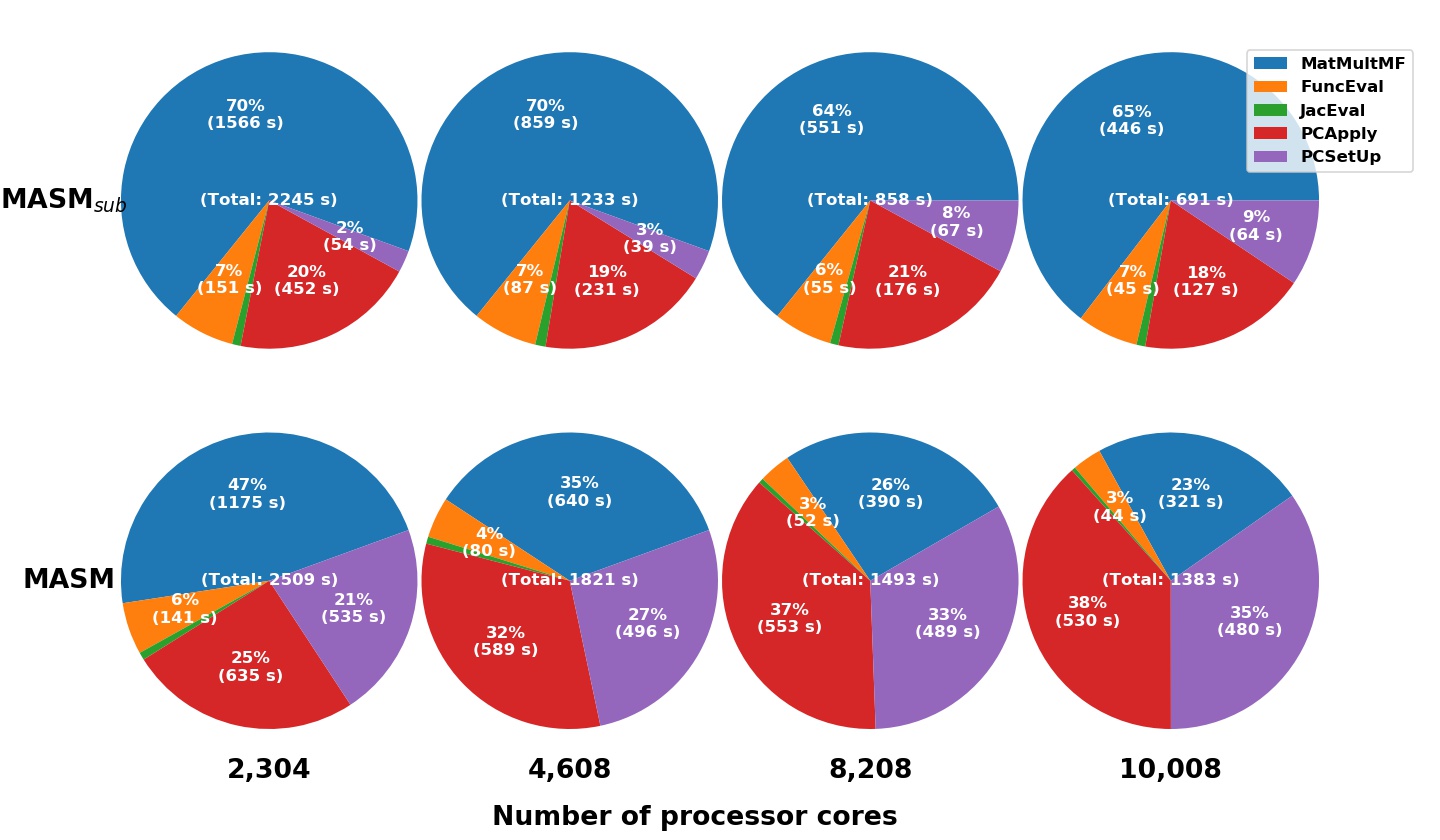} 
  \caption{Performance comparison with MASM for the problem with 1,253,400,480 unknowns. The figure shows the compute time spent on Jacobian evaluations (``JacEval"), function evaluations (``FuncEval"), matrix vector multiplications   via matrix-free (``MatMultMF"),  preconditioner setup (``PCSetUp") and preconditioner apply (``PCApply"), respectively.   \label{fig:pieplot_16_s4}}
\end{figure}

\begin{figure}
\centering
\includegraphics[width=0.32\linewidth]{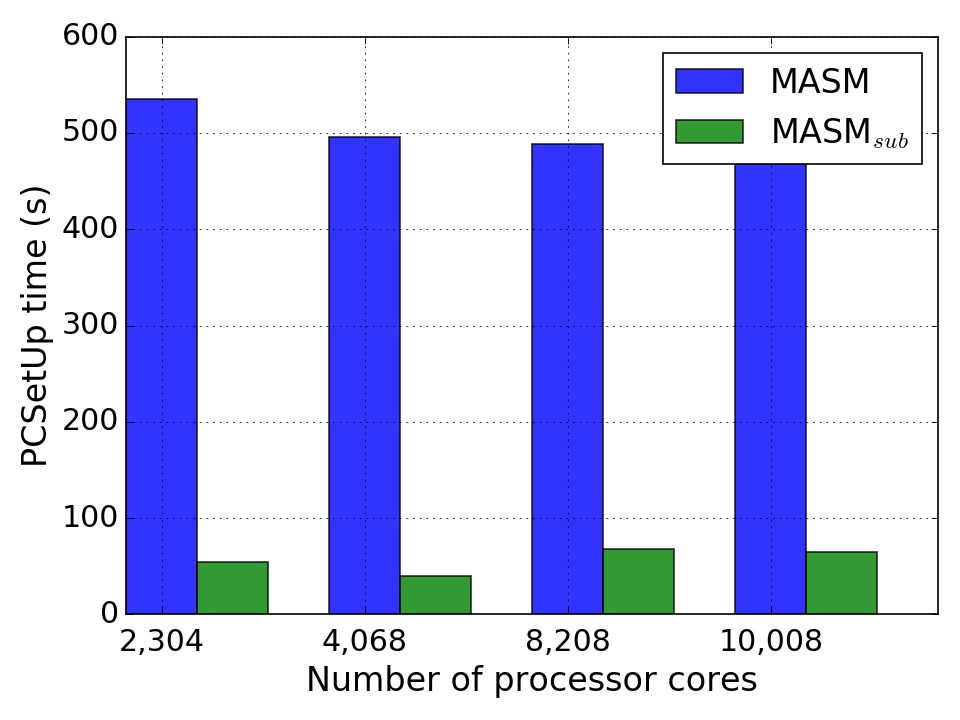} 
\includegraphics[width=0.32\linewidth]{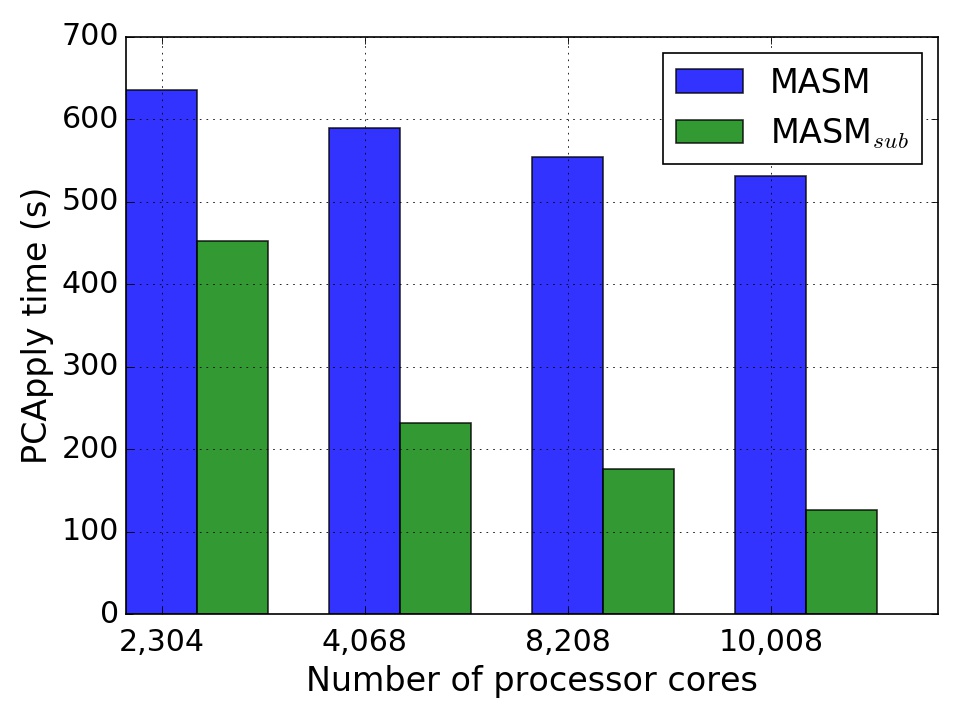} 
\includegraphics[width=0.32\linewidth]{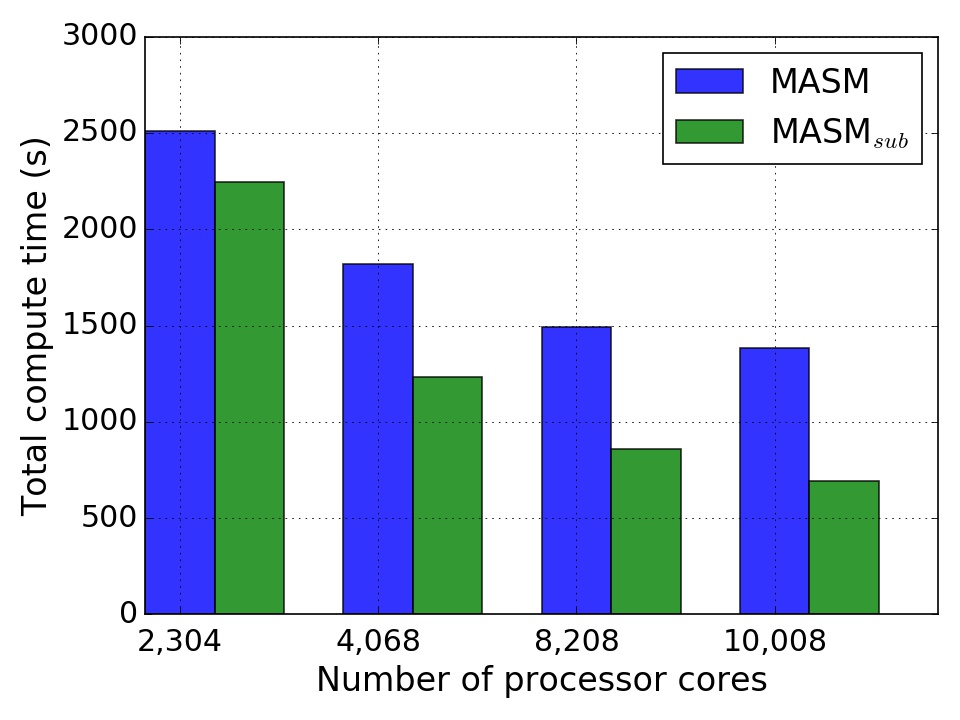}
\caption{The compute time comparison on different phases for the problem with 1,253,400,480 unknowns.  Left: preconditioner setup time; middle: preconditioner apply time; right: total compute time. \label{fig:barplot_16_s4}}
\end{figure}

To further verify the effectiveness of the proposed algorithm, we add more angular directions for each mesh node, which leads to 288 variables (24 angular directions $\times$ 12 energy groups) one   each  mesh node. The same mesh as before is used, but the resulting system is much larger than the previous test, having  1,253,400,480 unknowns, since more angular  directions are added for each mesh node.  The numerical results are summarized in Table~\ref{tab:amg_sub_coarse_s4}, Fig.~\ref{fig:pieplot_16_s4} and~\ref{fig:barplot_16_s4}.  In Table~\ref{tab:amg_sub_coarse_s4}, 
it is easily seen that the number of Newton iterations stays close to a constant as we increase the number of processor cores from 2,304 to 10,008, and the number of GMRES iterations also keeps as
a constant for both MASM and MASM$_{\text{sub}}$.  The number of GMRES iterations for  MASM$_{\text{sub}}$ is more than that obtained via MASM, but the performance of  MASM$_{\text{sub}}$ is not affected much and it is still much better than MASM.  The overall algorithm based on MASM$_{\text{sub}}$ has the parallel efficiency of $76\%$ at 10,008 cores while that for MASM is only $37\%$.  Similarly, from Fig.~\ref{fig:pieplot_16_s4} and~\ref{fig:barplot_16_s4},  we observed that the compute time of the preconditioner for MASM$_{\text{sub}}$ is much smaller than that for MASM.
\begin{table}
\scriptsize
\centering
\caption{Performance comparison with MASM for a problem with 1,253,400,480 unknowns. The resulting system of nonlinear equations with 1,253,400,480 unknowns 
is solved by an inexact Newton with MASM$_{\text{sub}}$ and  MASM  on  2,304, 4,608,  8,208, and 10,008 processor cores, respectively. \label{tab:amg_sub_coarse_s4}}
\begin{tabular}{c c c c c c  c c c c}
\toprule
$np$  &scheme& NI& LI& Newton& LSolver & MF   &  PCSetup & PCApply   & EFF \\
\midrule
2,304 & MASM$_{\text{sub}}$ &  13 & 191 &  2202 &  2027&1566& 54 &452 & 100\% \\
2,304 & MASM &  12 & 147 &  2466 & 2302&1176& 535 &  635 & -- \\
\midrule
4,608 &  MASM$_{\text{sub}}$  & 13 & 183& 1199& 1096 & 860 &41 & 232 & 92\%\\\
4,608 &  MASM  & 12 & 139&  1791& 1694 & 641 &496 & 589 & 61\%\\
\midrule
8,208 &  MASM$_{\text{sub}}$  & 13 & 183& 828& 763 &  552 &68 & 176 & 75\%\\
8,208 &  MASM  & 12 & 135& 1474& 1412 & 391 & 489 &  554 & 42\%\\
\midrule
10,008 &  MASM$_{\text{sub}}$  & 14 & 184& 672& 617 & 447 &65 & 127 & 76\%\\
10,008 &  MASM  & 12 & 134& 1369& 1317 & 322 &480 &  531 & 37\%\\
\bottomrule
\end{tabular}
\end{table}
%% Pie plot
%% bar plots

\subsection{Node balance improvement}
\begin{figure}
\centering
\includegraphics[width=0.32\linewidth]{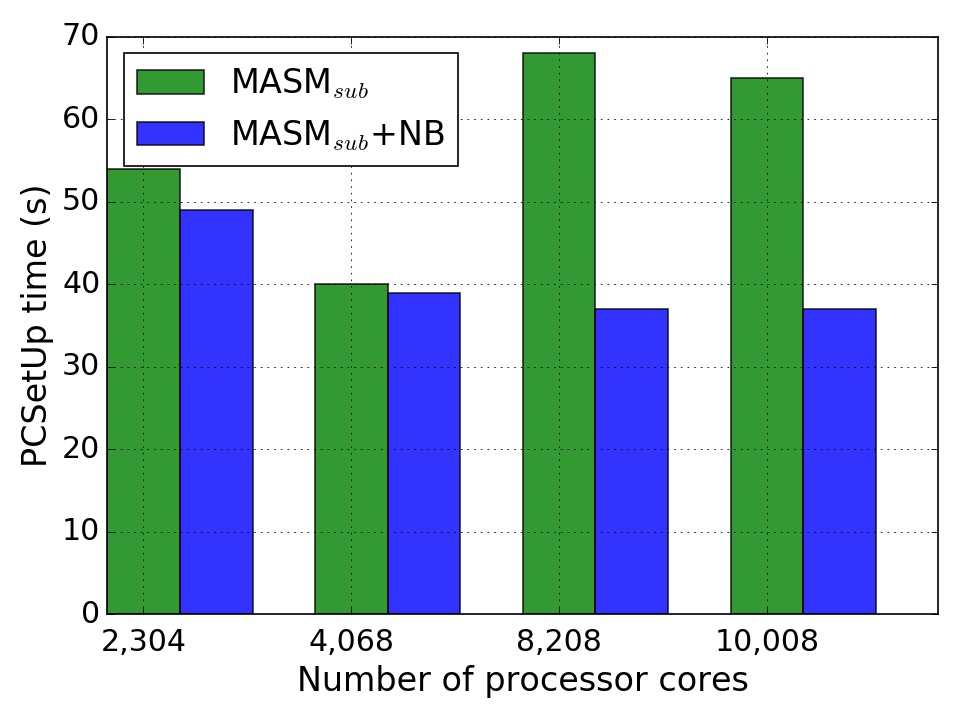} 
\includegraphics[width=0.32\linewidth]{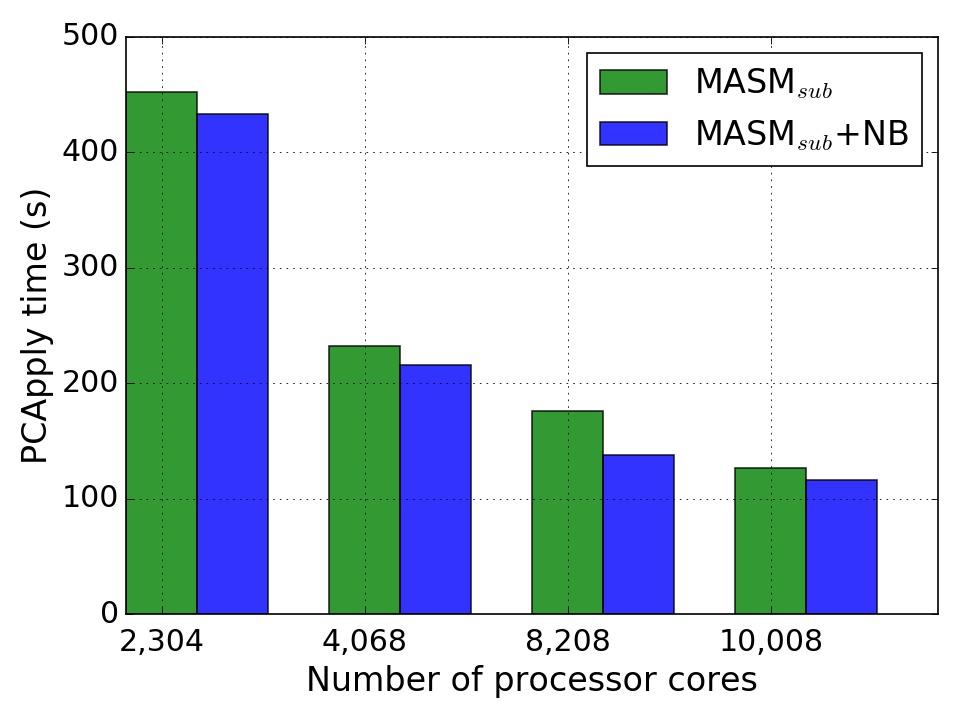} 
\includegraphics[width=0.32\linewidth]{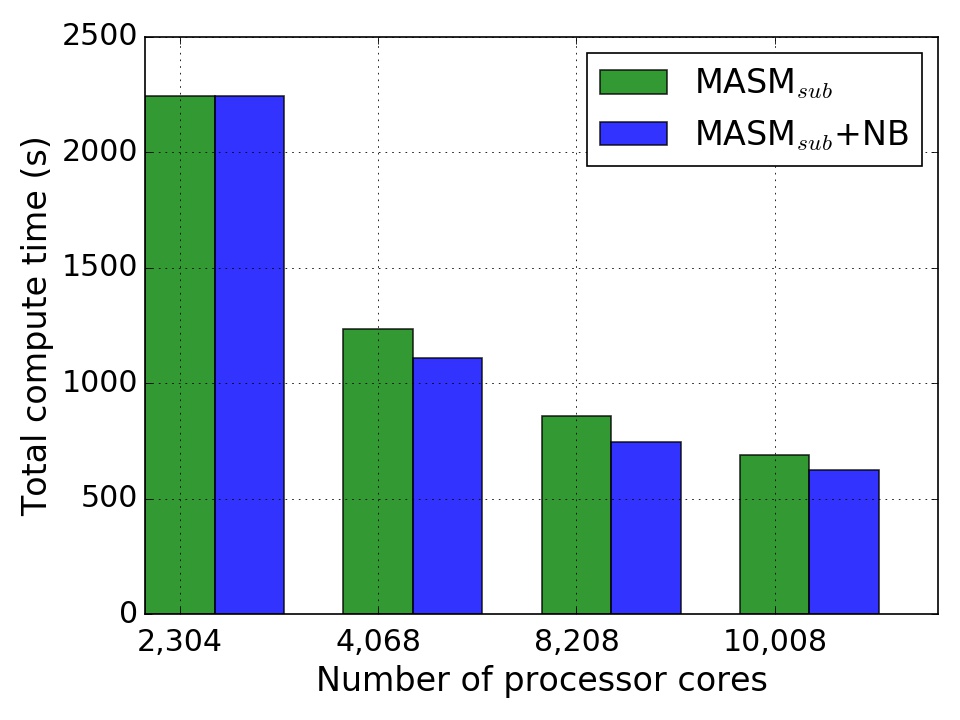}
\caption{The compute time comparison on different phases for the problem with 1,253,400,480 unknowns using a node balancing strategy.  Left: preconditioner setup time; middle: preconditioner apply time; right: total compute time. \label{fig:barplot_16_s4_nd}}
\end{figure}

Typically,  before a finite element simulation starts,  a dual graph of mesh (where each graph vertex corresponds to a mesh element) is partitioned into $np$ submeshes that have nearly  equal number of elements.  While the number of elements assigned to each processor core is nearly equivalent, some processor cores may have more mesh nodes  since the shared mesh nodes along the processor boundaries  are simply assigned to the cores  with lower MPI rank.  A scalable calculation requires  to balance both mesh elements and mesh nodes.  We use a partition-based node assignment, as discussed in our previous work \cite{kong2018general},  to balance the overall calculation.  The basic idea of the partition-based node assignment is that the processor boundary mesh is partitioned 
into two parts, and each part is assigned to a processor core who shares the processor boundary mesh with the other processor core. Interested readers are referred to \cite{kong2018general} for more details. The same configuration as before is used, and the numerical results are shown in Table~\ref{tab:nodeassignment_coarse_s2} and~\ref{tab:nodeassignment_coarse_s4}, and Fig.~\ref{fig:barplot_16_s4_nd}. 
\begin{table}
\scriptsize
\centering
\caption{Mesh node assignment for the problems with 417,800,160 unknowns. \label{tab:nodeassignment_coarse_s2}}
\begin{tabular}{c c c c c c  c c c c}
\toprule
$np$  &scheme& NI& LI& Newton& LSolver & MF   &  PCSetup & PCApply   & EFF \\
\midrule
1,152 & MASM$_{\text{sub}}$ &  13 & 191 & 1855 & 1701&1418& 26 &290 & -- \\
1,152 & MASM$_{\text{sub}}$+NB &  13 & 199 & 1821 & 1674&1407& 26 & 270 & 100\% \\
\midrule
2,304 &  MASM$_{\text{sub}}$  & 13 & 193& 989& 908 & 749 &21 & 154 & 92\%\\
2,304 &  MASM$_{\text{sub}}$+NB  & 13 & 207& 1005& 928 & 769 &19 & 158 & 91\%\\
\midrule
4,608 &  MASM$_{\text{sub}}$  & 13 & 202& 581& 535 & 440 &18 & 90 & 78\%\\
4,608 &  MASM$_{\text{sub}}$+NB  & 13 & 216& 579& 534 & 436 &18 &  89 & 79\%\\
\midrule
8,208 &  MASM$_{\text{sub}}$  & 14 & 216& 404& 372 & 294 &22 & 66 & 63\%\\
8,208 &  MASM$_{\text{sub}}$+NB  & 14 & 217& 379& 348 & 272 &21 & 62 & 67\%\\
\bottomrule
\end{tabular}
\end{table}
From Table~\ref{tab:nodeassignment_coarse_s2}, we observed that the compute time for different components  is further improved using the node assignment strategy.  The improvement becomes more obvious  when the number of variables for each mesh node is increased.  In Table~\ref{tab:nodeassignment_coarse_s4}, the preconditioner setup time is significantly reduced, for example, it is reduced to $37~s$
from $65~s$ at 10,008 processor cores, which leads to the parallel efficiency  increased from $76\%$ from $83\%$.  The same observation is found in Fig.~\ref{fig:barplot_16_s4_nd} as well.
\begin{table}
\scriptsize
\centering
\caption{Node assignment for the problem with with 1,253,400,480 unknowns.\label{tab:nodeassignment_coarse_s4}}
\begin{tabular}{c c c c c c  c c c c}
\toprule
$np$  &scheme& NI& LI& Newton& LSolver & MF   &  PCSetup & PCApply   & EFF \\
\midrule
2,304 &  MASM$_{\text{sub}}$ &  13 & 191 &  2202 &  2027 &1566& 54&452 & 100\% \\
2,304 &MASM$_{\text{sub}}$+ NB &  13 & 191 &  2203 & 2018&1575& 49&  433 & -- \\
\midrule
4,608 &   MASM$_{\text{sub}}$  & 13 & 183& 1199& 1096 & 860 &41 & 232 & 92\%\\
4,608 &  MASM$_{\text{sub}}$+NB  & 13 & 185&  1093& 1001 & 766 &39 & 216 & 100\%\\
\midrule
8,208 &   MASM$_{\text{sub}}$  & 13 & 183& 828& 763 &  552 &68 & 176 & 75\%\\
8,208 &  MASM$_{\text{sub}}$+NB  & 13 & 189& 732& 670 & 511 & 37 &   138 & 84\%\\
\midrule
10,008 &  MASM$_{\text{sub}}$  & 13 & 184& 672& 617 & 447 &65 & 127 & 76\%\\
10,008 &  MASM$_{\text{sub}}$+NB  & 13 & 187& 610& 557 & 420 &37 &  116 & 83\%\\
\bottomrule
\end{tabular}
\end{table}
%% bar plots

\subsection{Aggressive coarsening}
The numerical results shown earlier are obtained using 10 levels of aggressive coarsening. The aggressive coarsening, Alg.~\ref{alg:AHCIJP_subspace}, is used to reduce the complexities of the operators and the interpolations. More levels are applied by the aggressive coarsening, and the complexities  of the operators and the interpolations become lower, but at the same time the convergence may be deteriorated.  In this test, different numbers of  levels  of  aggressive coarsening  are applied in  MASM$_{\text{sub}}$, and the results are summarized in Table~\ref{tab:aggressive_coarsening}, where   ``Agg" denotes the number of levels of aggressive coarsening, and  ``Agg=0" corresponds to no aggressive coarsening.  
\begin{table}
\scriptsize
\centering
\caption{Different number of aggressive coarsening levels for  the problems with 417,800,160 unknowns. ``Agg" denotes the number of aggressive coarsening levels. \label{tab:aggressive_coarsening}}
\begin{tabular}{c c c c c c  c c c c}
\toprule
$np$  &Agg& NI& LI& Newton& LSolver & MF   &  PCSetup & PCApply   & EFF \\
\midrule
1,152 & 0  &  13 & 209  &  2242  &  2088 &  1853& 55&509 & 100\% \\
1,152 & 2 &  14 & 210 &  2068 & 1906    &1570&33 &347 & -- \\
1,152 & 4 &  13 & 206 &  2007 & 1853    & 1542&26 &322 & -- \\
1,152 & 10 &  13 & 191 & 1855 & 1701   &1418& 26 &290 & -- \\
\midrule
2,304 &  0  & 13 & 208& 1199&1118 & 820 &45 &  278 & 77\%\\
2,304 &  2  & 13 & 195& 1027& 945 & 758 &25 & 178 & 90\%\\
2,304 &  4  & 13 & 206& 1059&978 & 801 &20 & 174 & 88\%\\
2,304 &  10  & 13 & 193& 989& 908 & 749 &21 & 154 & 94\%\\
\midrule
4,608 &  0  & 13 & 206& 677&632 &  457 & 34 &  160 & 68\%\\
4,608 &  2  & 13 & 199& 591& 546 & 423 &21 & 113 & 78\%\\
4,608 &  4  & 13 & 210&  606&558 & 449 &18 & 104 & 80\%\\
4,608 &  10  & 13 & 202& 581& 534 & 440 &18 & 90 & 94\%\\
\midrule
8,208 &  0  & 14 & 217& 490&457 &  302 & 46 &  123 & 53\%\\
8,208 &  2  & 13 & 204& 412& 379 & 277 &25 & 87 & 63\%\\
8,208 &  4  & 13 & 200&  385& 355 & 273 &21 & 70 & 68\%\\
8,208 &  10  & 14 & 216& 404& 372 & 294 &22 & 66 & 64\%\\
\bottomrule
\end{tabular}
\end{table}
From Table~\ref{tab:aggressive_coarsening}, we observed that both Newton iteration and GMRES iteration  do not change much when we apply different numbers of levels  of the aggressive coarsening  in MASM$_{\text{sub}}$. The preconditioner time including the setup and apply is reduced significantly when the number  of levels of of aggressive coarsening  is increased from 0 to 2, and then it slightly decreases  when we increase ``Agg" from 2, 4 to 10. For example, in 8,208-core case, the preconditioner apply time is decreased from $123~s$  to $87~s$ when ``Agg" is increased from 0 to 2, and it 
is reduced to $70~s$ and $66~s$ when ``Agg" is 4 and 8.  Due to these factors, the corresponding parallel efficiency is increased from $53\%$ to $64\%$ when the number of levels of aggressive coarsening  is  increased from 0 to 10. It is the reason why we have used  10 levels of aggressive coarsening  in our  previous tests. 

\subsection{Strong scalability}
\begin{figure}
\centering
\includegraphics[width=0.45\linewidth]{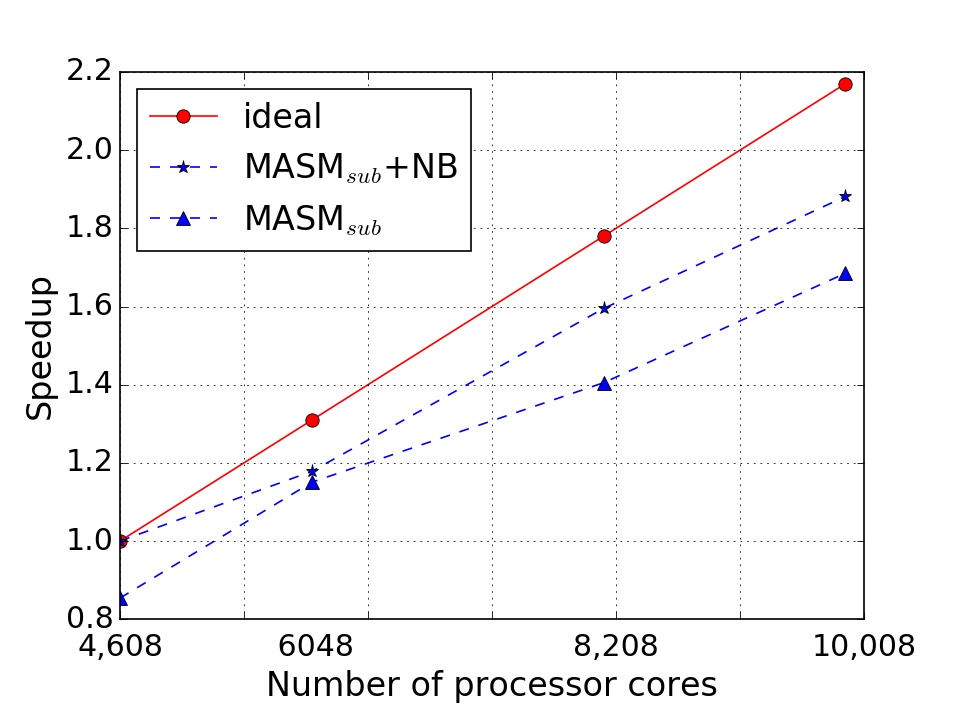} 
\includegraphics[width=0.45\linewidth]{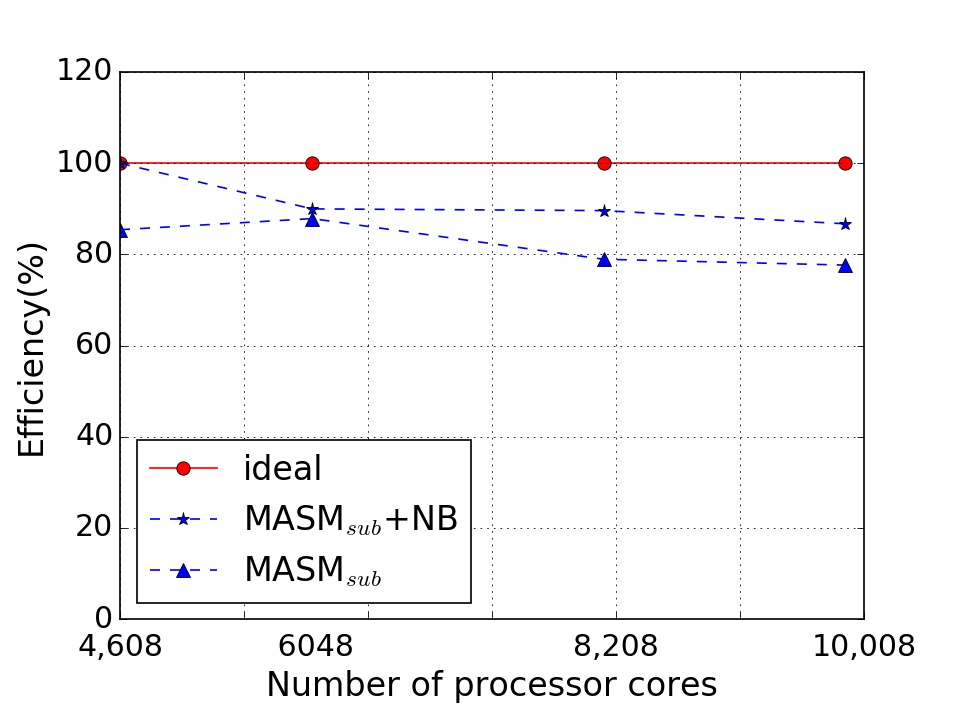} 
\caption{Parallel efficiency and speedup for the problem with 2,482,224,480 unknowns. Left: speedup, right: parallel efficiency.  \label{fig:parallel_efficiency_100_s2_nd}}
\end{figure}

In this test, we study the strong scalability using a ``fine"  mesh with 25,856,505 nodes and 26,298,300 elements. The resulting eigenvalue system of equations  with 2,482,224,480 unknowns is solved by an inexact Newton preconditioned by MASM$_{\text{sub}}$.  At the beginning  of the strong scaling study, we also test the impact of the number of levels of aggressive coarsening  on the overall algorithm for the ``fine" mesh case.  4 and 10 levels of aggressive coarsening  are tested, and the results are summarized in Table~\ref{tab:aggressive_coarseing_100_s2}.  The case with 10 levels of aggressive coarsening  has slightly better results  than  that obtained with ``Agg=4".  
\begin{table}
\scriptsize
\centering
\caption{Different number of aggressive coarsening levels for the problem   with 2,482,224,480 unknowns. \label{tab:aggressive_coarseing_100_s2}}
\begin{tabular}{c c c c c c  c c c c}
\toprule
$np$  &Agg& NI& LI& Newton& LSolver & MF   &  PCSetup & PCApply   & EFF \\
\midrule
4,608 & 4 &  12 & 161  &  2596  &  2369 &  1853& 114&457 & 100\% \\
4,608 & 10 &  13 & 171 &  2808 & 2567&2002&131 &511 & -- \\
\midrule
6,048 &  4  & 13 & 170& 2125&1939 & 1499 &95 &  391 & 93\%\\
6,048 &  10  &13& 171& 2084&1898&  1510&  84 & 352& 95\%\\
\midrule
8,208 & 4  & 12&160 &1749& 1613& 1134& 150&404 & 83\%\\
8,208 & 10  & 12&161& 1707& 1570 &1139 & 142&  364& 85\%\\
\midrule
10,008 & 4  &12& 161&  1459 &  1343 &949&141& 308& 82\%\\
10,008 & 10  & 12 & 160& 1423 & 1307&   951& 138& 274& 84\%\\
\bottomrule
\end{tabular}
\end{table}
We therefore use 10 levels of aggressive coarsening in the following scaling study.  The numerical results are shown in Table~\ref{tab:strong_scaling_100_s2}. The performance of the node balancing  strategy  is also reported, and the corresponding algorithm is denoted  as ``MASM$_{\text{sub}}$+NB".  Again, for MASM$_{\text{sub}}$, the numbers of Newton iterations and GMRES iterations stay as constants when we increase the number of processor cores, which indicates  that the algorithm is mathematically scalable.  The total compute time (``Newton") is decreased proportionally   when 
we increase the number of processor cores from 4,608 to 10,008. For instance, the total compute is reduced from $2808~s$ to $2084~s$,  when we increase the number of processor cores from 
4,608 to 6,048, and it further is reduced to $1707~s$ and $1423~s$  when the number of processor cores is 8,208 and 10,008.  The preconditoner setup is not scalable, but it does not affect the overall performance much since it accounts for only a small portion of the total compute time. A good parallel efficiency of $77\%$ is achieved at 10,008.  While the performance of MASM$_{\text{sub}}$ is already good, it can be further enhanced  by applying the node balancing  strategy to make the overall distribution more balanced.  In the odd rows of Table~\ref{tab:strong_scaling_100_s2},
the node balancing  strategy is able to improve the overall algorithm performance, especially, the preconditioner setup time is significantly  reduced.  For example, at 8,208 cores, the compute time is reduced by $200~s$ when the node balancing  strategy is used, most of the time reduction results from the improvement of the preconditioner setup. At 10,008 cores, the preconditioner setup time is 
reduced to $82~s$ from $138~s$.  The preconditioner apply  also benefits from a more balanced workload, for example, the preconditioner apply time is reduced from $274~s$ to $208~s$ at 10,008 
processor cores.  We have almost-perfect  parallel efficiencies  when we use MASM$_{\text{sub}}$+NB.  The parallel efficiency  is as high as $87\%$, even when the number of processor cores is more than 10,000.   The Schwarz preconditioner together with both the subspace-based coarsening and the partition-based  node balancing   offers a highly scalable solver for the eigenvalue calculations  for the targeting application. The corresponding  parallel efficiency and speedup are plotted in Fig.~\ref{fig:parallel_efficiency_100_s2_nd}. 
\begin{table}
\scriptsize
\centering
\caption{ Strong scalability  with a ``fine" mesh  with 25,856,505 nodes,  26,298,300 elements, and 2,482,224,480 unknowns.  \label{tab:strong_scaling_100_s2}}
\begin{tabular}{c c c c c c  c c c c c}
\toprule
$np$  &Scheme& NI& LI& Newton& LSolver & MF   &  PCSetup & PCApply   &NR& EFF \\
\midrule
4,608 & MASM$_{\text{sub}}$+NB &  12 & 160  &  2398  &  2182 &  1743& 94&396&1.8& 100\% \\
4,608 & MASM$_{\text{sub}}$ &  13 & 171 &  2808 & 2567&2002&131 &511 &2.2& -- \\
\midrule
6,048 &  MASM$_{\text{sub}}$+NB  & 13 & 176& 2035&1858 & 1470 &87 &   337 &1.8& 90\%\\
6,048 &  MASM$_{\text{sub}}$  &13& 171& 2084&1898&  1510&  84 & 352& 2.7&89\%\\
\midrule
8,208 & MASM$_{\text{sub}}$+NB  & 12&167 &1504& 1373& 1081& 76&244 &2& 90\%\\
8,208 & MASM$_{sub}$  & 12&161& 1707& 1570 &1139 & 142&  364&2.8& 79\%\\
\midrule
10,008 & MASM$_{\text{sub}}$+NB  &13& 168&  1275 &   1160 &896&82& 208&2& 87\%\\
10,008 & MASM$_{sub}$  & 12 & 160& 1423 & 1307&   951& 138& 274&2.9&77\%\\
\bottomrule
\end{tabular}
\end{table}
%% bar plots

\section{Conclusions}
A parallel Newton-Kyrlov-Schwarz method is studied for the numerical simulation of the multigroup  neutron transport equations on 3D unstructured spatial meshes. A hierarchal partitioning is used to divide the computational domain into a large number of subdomains while the existing  partitioners are from ideal. Two novel components including the subspace-based coarsening and the partition-based workload balancing are introduced and carefully studied. The total compute time using the subspace-based  coarsening is halved   compared with the traditional approach, and therefore the corresponding parallel efficiency is doubled    to $76\%$ when more than 10,000 processor cores are used. In addition, the partition-based workload balancing is able to assign a nearly  equal amount of work to each processor core, and the parallel efficiency at 10,000 is further increased  to $87\%$.  The scalability of the overall algorithm is studied for  a realistic application with 2,482,224,480 unknowns on a supercomputer with up to 10,008 processor cores.

While this paper focuses on the multilevel Schwarz preconditioner, other popular methods such as DSA \cite{alcouffe1977diffusion} and  NDA (nonlinear diffusion acceleration method) \cite{schunert2017flexible} will be explored in the future.
\bibliographystyle{siamplain}
\bibliography{hmg}

\end{document}